\documentclass{article}
\usepackage{amssymb}
\usepackage{eurosym}
\usepackage{epsfig}
\usepackage{amsmath}
\usepackage{amsfonts}
\usepackage{pgf,tikz}
\usepackage{enumerate}
\usepackage{tabularx}
\usepackage{array}
\usepackage{booktabs}
\usepackage{microtype}
\usepackage{mathrsfs}
\usepackage{cite}
\usepackage{subfigure}
\usepackage[colorlinks=true, linkcolor=black, citecolor=black, urlcolor=black]{hyperref}
\usepackage[top=2cm, bottom=2.5cm, left=2.5cm, right=2.5cm, headsep=1cm]{geometry}

\newtheorem{theorem}{\sc Theorem}[section]
\newtheorem{proposition}{\sc Proposition}[section]
\newtheorem{lemma}{\sc Lemma}[section]

\newtheorem{remark}{\sc Remark}[section]
\newtheorem{corollary}{\sc Corollary}[section]
\newtheorem{example}{\sc Example}[section]

\def\tanh#1{\,{\normalsize tanh}{\,#1}\,}

\def\qed{\hbox to 0pt{}\hfill$\rlap{$\sqcap$}\sqcup$\medbreak}

\title{Positive coexistence solutions for systems of multi-kernel Hammerstein type equations}

\author{Laura Mª Fernández-Pardo$\,^{1,2}$ and Jorge Rodríguez-López$\,^{1,2}$} 

\date{}
\begin{document}
 \maketitle

\begin{center}  {\small $^1$ CITMAga, 15782, Santiago de Compostela, Spain.  \\
		 $^2$ Departamento de Estat\'{\i}stica, An\'alise Matem\'atica e Optimizaci\'on, \\ Universidade de Santiago de Compostela, \\ 15782, Facultade de Matem\'aticas, Campus Vida, Santiago, Spain.\\  Email: laura.fernandez.pardo@usc.es; jorgerodriguez.lopez@usc.es}
\end{center}

\medbreak

\noindent {\it Abstract.} We investigate positive solutions with both components nontrivial for systems of Hammerstein integral equations associated with a general $2\times2$ matrix of Green's kernels, thereby extending the decoupled (diagonal) case to systems with cross-coupling. Such systems arise from boundary value problems in which the unknowns may be coupled both in the differential equations and in the boundary conditions, a setting that has received considerably less attention than the decoupled case. Using a component-wise Krasnosel'ski\u{\i} compression--expansion fixed point theorem in product cones, we derive existence and non-existence criteria. A multiplicity result guaranteeing at least two positive solutions is also obtained, by means of the fixed point index in cones. As an application, we consider radially symmetric solutions for an elliptic system on an annulus with coupled Robin boundary conditions.

\medbreak

\noindent     \textit{2020 MSC:} 45G15, 35J61, 34B08, 34B18, 47H10.

\medbreak

\noindent     \textit{Key words and phrases.} Coupled boundary value systems; positive solution; nonlinear systems; Hammerstein integral equations; radial elliptic systems.

\section{Introduction}
Positive solutions of nonlinear boundary value problems have been the subject of intensive research for several decades, both because of their intrinsic mathematical interest and because of their role in modeling steady states in fluid mechanics, chemical reactor theory and population dynamics (see, e.g., \cite{amann, Deim, guolak,Leung}). A standard approach to such problems is to use Green's functions to formulate them as Hammerstein integral equations. In this paper, we study the existence, non-existence, and multiplicity of positive solutions, with both components nontrivial, for the Hammerstein system of the form
\begin{equation}
	\label{1.1}
	\left\{
	\begin{aligned}
		u(t) &= \lambda_1\int_{I} G_1^1(t,s)f_1(s,u(s),v(s))\,ds
		+ \lambda_2\int_{I} G_2^1(t,s)f_2(s,u(s),v(s))\,ds,
		&& t\in I,\\[1em]
		v(t) &= \lambda_1\int_{I} G_1^2(t,s)f_1(s,u(s),v(s))\,ds
		+ \lambda_2 \int_{I} G_2^2(t,s)f_2(s,u(s),v(s))\,ds,
		&& t\in I,
	\end{aligned}
	\right.
\end{equation}
where $\lambda_i>0$, $I:=[A,B]$ with $A<B$, and the kernels $G^j_i:I\times I\to\mathbb{R}$ and the nonlinearities $f_i:I\times[0,\infty)^2\to[0,\infty)$, $i,j=1,2$, are continuous. 

Systems of the form \eqref{1.1} arise from coupled systems of two differential equations, of orders $n,m\in\mathbb{N}$, depending on both $u$ and $v$, with nonlinear right-hand sides $f_i$ and subject to $n+m$ nonlocal linear boundary conditions. More precisely, fixed points of \eqref{1.1} correspond to solutions of the boundary value problem
\begin{equation}
	\label{1.2}
	\left\{
	\begin{aligned}
		\mathcal{L}_1(u,v)(t) &= \lambda_1 f_1(t,u(t), v(t)), \quad t \in I, \\
		\mathcal{L}_2(u,v)(t) &= \lambda_2 f_2(t, u(t), v(t)), \quad t \in I, \\
		\mathcal{B}_i(u,v) &= \delta_i \mathcal{C}_i(u,v), \quad i = 1, \dots, n+m,
	\end{aligned}
	\right.
\end{equation}
where the linear differential operators $\mathcal{L}_1$ and $\mathcal{L}_2$ are defined by
\begin{equation}
	\begin{aligned}
		\mathcal{L}_1(u,v)(t) &:= u^{(n)}(t) + \sum_{k=0}^{n-1} a_k u^{(k)}(t) + \sum_{i=0}^{m-1} b_i v^{(i)}(t), \\
		\mathcal{L}_2(u,v)(t) &:= v^{(m)}(t) + \sum_{k=0}^{n-1} c_k u^{(k)}(t) + \sum_{i=0}^{m-1} d_i v^{(i)}(t).
	\end{aligned}
\end{equation}
Here, $\delta_i \in \mathbb{R}$ for all $i=1, \dots, n+m$, while $a_k, c_k \in \mathbb{R}$ for $k=0, \dots, n-1$ and $b_i, d_i \in \mathbb{R}$ for $i=0, \dots, m-1$. Furthermore, $\mathcal{C}_i : C^n(I) \times C^m(I) \to \mathbb{R}$ are linear continuous operators, and $\mathcal{B}_i$ are general two-point linear boundary operators given by
\[
\mathcal{B}_i(u,v) = \sum_{j=0}^{n-1} \left( \alpha_j^i u^{(j)}(A) + \beta_j^i u^{(j)}(B) \right) + \sum_{l=0}^{m-1} \left( \tilde{\alpha}_l^i v^{(l)}(A) + \tilde{\beta}_l^i v^{(l)}(B) \right), \quad i=1, \dots, n+m,
\]
where $\alpha_j^i, \beta_j^i, \tilde{\alpha}_l^i, \tilde{\beta}_l^i \in \mathbb{R}$ for $j=0, \dots, n-1$, $l=0, \dots, m-1$, and $i=1, \dots, n+m$. The equivalence between \eqref{1.2} and \eqref{1.1} was recently established in \cite{Mou}, where the authors rewrote \eqref{1.2} as a first-order matrix system and determined explicitly the associated Green's matrix. Concrete systems of boundary value problems with coupled boundary conditions had been previously considered, for instance, in \cite{AK,CS,IMP}.

When neither the differential operators nor the boundary operators mix the unknowns, that is, when $\mathcal{L}_1$ depends only on $u$, $\mathcal{L}_2$ depends only on $v$, and each $\mathcal{B}_i$ involves only one of the two variables, the off-diagonal kernels vanish and \eqref{1.1} reduces to the diagonal system
\begin{equation}
	\label{1.4}
	\left\{
	\begin{aligned}
		u(t) &= \lambda_1 \int_{I} G_1^1(t,s)f_1(s,u(s),v(s))\,ds,
		&& t\in I,\\[1em]
		v(t) &= \lambda_2 \int_{I} G_2^2(t,s)f_2(s,u(s),v(s))\,ds,
		&& t\in I.
	\end{aligned}
	\right.
\end{equation}
The vast majority of the literature on boundary value problems (see, e.g., \cite{sousa, LFP_JRL, PrecupFPT,JRL, IP,Y} and the references therein) is concerned with \eqref{1.4}: each component is governed by its own scalar kernel, and the interaction between $u$ and $v$ takes place only through the nonlinearities $f_1, f_2$. In contrast, when coupling is also present at the level of the derivatives, and in particular in the boundary conditions, the integral formulation involves the full Green's matrix
	\begin{equation*}
		\label{eq:intro_green_matrix}
		\mathbb{G}(t,s)
		=
		\begin{pmatrix}
			G_1^1(t,s) & \,\,G_2^1(t,s) \\[0.6em]
			G_1^2(t,s) & \,\,G_2^2(t,s)
		\end{pmatrix},
	\end{equation*}
	whose four entries are, in general, nonzero, non-symmetric, and structurally different from one another. Boundary conditions of this type arise naturally in applications where the state of one component at the boundary is influenced by the other, as in boundary interactions between species in reaction--diffusion models (see, for instance, the prey--predator system with coupled boundary conditions studied by Leung \cite{Leung}). Nevertheless, as far as the authors are aware, a systematic study of the existence of positive solutions for such systems has only been addressed by means of the abstract approach developed in \cite{Mou}.

A problem exemplifying this situation, which will be used throughout the paper to illustrate our abstract results, is the elliptic system posed on the annulus $\Omega = \{x \in \mathbb{R}^3 :  A < |x| < B\}$, with $0 < A<B<+\infty$, $B-A=1$,
	\begin{equation}\label{eq:intro_pde_system}
		\left\{
		\begin{aligned}
			-\Delta U &= \lambda_1 h_1(|x|)
			f_1\big(U, V\big), & x \in \Omega, \\[0.6em]
			-\Delta V &= \lambda_2 h_2(|x|)
			f_2\big(U, V\big), & x \in \Omega,
		\end{aligned}
		\right.
	\end{equation}
	subject to the coupled Robin boundary conditions
	\begin{equation}\label{eq:intro_pde_bc}
		\left\{
		\begin{aligned}
			U(x) - \frac{A}{B} \frac{\partial U}{\partial \nu}(x)
			+ \alpha V(x) &= 0, \quad |x| = A; \qquad &
			U(x) + \frac{B}{A} \frac{\partial U}{\partial \nu}(x)
			&= 0, \quad |x| = B, \\[0.8em]
			V(x) - \beta \frac{A}{B} \frac{\partial V}{\partial \nu}(x)
			+ \alpha U(x) &= 0, \quad |x| = A; \qquad &
			V(x) + \beta \frac{B}{A} \frac{\partial V}{\partial \nu}(x)
			&= 0, \quad |x| = B,
		\end{aligned}
		\right.
	\end{equation}
	where $\partial/\partial\nu$ denotes differentiation in the radial direction, $\alpha \in \mathbb{R}$ is a coupling parameter, and $\beta \ge 0$. Radial solutions of elliptic equations and systems on annuli have been widely investigated by topological methods; see, for instance, \cite{LanWebb,LanLin,CiaPie,wang,doOLU,doOLU2,G1,IP}. In these works, however, each boundary condition involves a single unknown, so that the associated Green's matrix is diagonal and the problem falls within the class described by \eqref{1.4}. In contrast, when $\alpha\neq0$, the inner boundary conditions in \eqref{eq:intro_pde_bc} couple $U$ and $V$, and the associated Green's matrix is no longer diagonal. Positive solutions of elliptic systems on annuli with coupled boundary conditions have received relatively little attention in the literature, with the notable exception of \cite{CIP}. Nevertheless, the fixed point formulation proposed in \cite{CIP} differs from the one employed here, since it relies on a system of perturbed Hammerstein type equations, instead of a system of the form \eqref{1.1}. In the appendices, we compute explicitly the Green's matrix of the radial problem associated with \eqref{eq:intro_pde_system}--\eqref{eq:intro_pde_bc}.

The presence of nonzero off-diagonal kernels in \eqref{1.1} introduces analytical difficulties absent in the decoupled case \eqref{1.4}. Indeed, since each component is governed by the sum $\lambda_1 G_1^i f_1 + \lambda_2 G_2^i f_2$, establishing existence requires balancing the simultaneous contribution of both nonlinearities in each equation (see Theorem \ref{th2.1}). As a consequence, their asymptotic growth conditions become coupled, leading to structural obstructions that make certain combinations of sublinear and superlinear growth incompatible (see Remark~\ref{rem:jde_compatibility}). Even so, we provide directional and diagonal asymptotic conditions at zero and at infinity which ensure the existence of at least one positive solution for \eqref{1.1}.

The main contributions of this work can be summarized as follows:
\begin{enumerate}[(a)]
	\item \textbf{Vectorial localization without semi-trivial solutions:}
	By employing the component-wise Krasnosel'ski\u{\i}--Guo fixed point theorem from \cite{LFP_JRL}, we establish existence results in which both components are strictly positive, thereby ruling out semi-trivial solutions of the form $(u,0)$ or $(0,v)$, unlike in \cite{Mou}, where such solutions may occur.
	
	\item \textbf{Generalization to coupled boundary value problems:}  We derive existence and non-existence results that address the coupling between the two equations, extending classical results for decoupled systems \eqref{1.4}, recovered when $G_2^1 \equiv G_1^2 \equiv 0$.
	
	\item \textbf{Multiplicity of two nontrivial solutions:} By combining superlinear growth at infinity with superlinear behavior along the diagonal near the origin, we establish a multiplicity theorem guaranteeing at least two distinct positive solutions with nontrivial components for sufficiently small values of $\lambda_1$ and $\lambda_2$. Roughly speaking, it provides a \textit{small} and a \textit{large} solution for system \eqref{1.1}. This multiplicity result is based on a novel computation of the fixed point index of the operator in the product space.
\end{enumerate}

The paper is organized as follows. In Section~2 we recall the vectorial compression--expansion Krasnosel'ski\u{\i}--Guo fixed point theorem and set up the operator formulation of \eqref{1.1}. Section~3 contains the existence and non-existence results, together with an application to \eqref{eq:intro_pde_system}--\eqref{eq:intro_pde_bc}. In Section~4 we develop the directional and diagonal asymptotic criteria at zero and infinity. Finally, Section~5 is devoted to a multiplicity result. The computation of the Green's matrix, of its extrema and of the Harnack constants is collected in Appendices~\ref{app:A}, \ref{app:B} and~\ref{app:C}.

\section{Preliminaries }
\label{sec:preliminaries}

In this section we fix the abstract setting of the paper and introduce the main ideas and concepts used throughout. We first recall the fixed point theorem on which most of our results rely, and then show how problem \eqref{1.1} fits into it.

Our main tool is the following component-wise version of the Krasnosel'ski\u{\i}--Guo theorem, established in \cite[Theorem 3.1]{LFP_JRL} (see also \cite{PrecupFPT} for the original vectorial version due to Precup).

\begin{theorem}
	\label{th1.1}
	Let $K_1$ and $K_2$ be two cones in the normed spaces $X_1$ and $X_2$, respectively.
	We denote $K := K_1 \times K_2$ the corresponding cone in
	$X = X_1 \times X_2$, and for $r,R \in \mathbb{R}^2$ with
	$0 < r_i < R_i$ $(i=1,2)$,
	\[
	\overline{K}_{r,R}
	:=
	\left\{
	x=(x_1,x_2)\in K :
	r_i \leq \|x_i\| \leq R_i
	\text{ for } i=1,2
	\right\}.
	\]
	
	Assume that $\mathcal{T}=(\mathcal{T}_1,\mathcal{T}_2):
	\overline{K}_{r,R} \to K$ is a compact map and for each
	$i\in\{1,2\}$ either of the following conditions holds in
	$\overline{K}_{r,R}$:
	\begin{enumerate}
		\item[(a)]
		$\|\mathcal{T}_i x\| \geq \|x_i\|$
		if $\|x_i\|=r_i$ and
		$\|\mathcal{T}_i x\| \leq \|x_i\|$
		if $\|x_i\|=R_i$, or
		
		\item[(b)]
		$\|\mathcal{T}_i x\| \leq \|x_i\|$
		if $\|x_i\|=r_i$ and
		$\|\mathcal{T}_i x\| \geq \|x_i\|$
		if $\|x_i\|=R_i$.
	\end{enumerate}
	Then $\mathcal{T}$ has at least one fixed point
	$x=(x_1,x_2)\in K$ such that $	r_i \leq \|x_i\| \leq R_i,$ $i=1,2$.
\end{theorem}

\begin{remark}
		\label{rem2.1}
In \cite{LFP_JRL}, it is further shown that, when the inequalities are strict, that is, when conditions $(a)$ and $(b)$ are replaced by, respectively,
\begin{enumerate}
	\item[$(a^*)$]$\|\mathcal{T}_i x\| > \|x_i\|$ if $\|x_i\|=r_i$ and $\|\mathcal{T}_i x\| < \|x_i\|$ if $\|x_i\|=R_i,$
	
	\item[$(b^*)$]$\|\mathcal{T}_i x\| < \|x_i\|$ if $\|x_i\|=r_i$ and $\|\mathcal{T}_i x\| > \|x_i\|$ if $\|x_i\|=R_i,$
\end{enumerate}
then not only does there exist at least one fixed point $x=(x_1,x_2)\in K$ of $\mathcal{T}$ satisfying $r_i < \|x_i\| < R_i,$ $ i=1,2$, but the fixed point index of $\mathcal{T}$ in $K_{r,R}:=\left\{x=(x_1,x_2)\in K : r_i < \|x_i\| < R_i,\ i=1,2\right\}$ is also determined explicitly by 
\[i_K(\mathcal{T},K_{r,R})=(-1)^s,\]
where $s\in\{0,1,2\}$ denotes the number of components of $\mathcal{T}$ satisfying condition $(b^*)$.

The knowledge of the fixed point index makes this result a powerful tool for establishing multiplicity results.
\end{remark}

In the sequel, following the usual terminology, it is said that the component  $\mathcal{T}_i$ of the operator $\mathcal{T}$ is \emph{compressive} if $(a)$ holds, whereas $\mathcal{T}_i$ is said to be \emph{expansive} in case that condition $(b)$ is satisfied. 

Theorem~\ref{th1.1} is particularly well suited to problem~\eqref{1.1}. Since the localization is given component by component, any fixed point in $\overline{K}_{r,R}$ satisfies $\|x_1\| \ge r_1 > 0$ and $\|x_2\| \ge r_2 > 0$, so semi-trivial solutions such as $(u,0)$ or $(0,v)$ are excluded. Moreover, the theorem allows compressive behaviour in one component to be combined with expansive behaviour in the other. To apply it, we take $X_1 = X_2 = C(I)$ with the supremum norm, and we need three ingredients: an operator $\mathcal{T}=(\mathcal{T}_1,\mathcal{T}_2)$ whose fixed points are the solutions of~\eqref{1.1}; cones $K_1,K_2$ that are invariant under $\mathcal{T}$; and a practical way of checking conditions $(a)$ or $(b)$ on the sets $\left\{x\in \overline{K}_{r,R}: \|x_i\| = r_i\right\}$ and $\left\{x\in \overline{K}_{r,R}: \|x_i\| = R_i\right\}$.

Accordingly, throughout this work the following assumptions will be imposed on the kernels and nonlinearities:
\begin{enumerate}
	\item[$(H_1)$] $G_j^i: I \times I \to \mathbb{R}^+ := [0, \infty)$ is continuous for all $i, j \in \{1, 2\}$.
	\item[$(H_2)$] For each $i \in \{1, 2\}$, there exist a subinterval $[a, b] \subseteq I$, constants $c_i \in (0, 1]$, $\tilde{c}_i \in [0, 1]$, and nonnegative functions $\Phi_j^i \in L^1(I)$ ($j = 1, 2$) such that:
	\begin{enumerate}
		\item[\emph{(i)}] For all $j \in \{1, 2\}$, the upper bounds hold:
		\begin{equation}\label{eq:H2_upper}
			G_j^i(t, s) \le \Phi_j^i(s), \quad \text{for all } (t, s) \in I \times I.
		\end{equation}
		
		\item[\emph{(ii)}] For all $j \in \{1, 2\}$, the lower inequalities hold:
		\begin{align}
			G_j^i(t, s) &\ge c_i \Phi_j^i(s), \quad \text{for all } t \in [a, b], \; s \in I, \label{eq:H2_lower_loc} \\[1ex]
			G_j^i(t, s) &\ge \tilde{c}_i \Phi_j^i(s), \quad \text{for all } t \in I, \; s \in I. \label{eq:H2_lower_glob}
		\end{align}
		
		\item[\emph{(iii)}] The diagonal kernel $G_i^i$ is strictly positive on $[a,b]\times [a,b]$.
	\end{enumerate}
	\item[$(H_3)$] $f_j: I \times \mathbb{R}^+ \times \mathbb{R}^+ \to \mathbb{R}^+$ is continuous for $j = 1, 2$.
\end{enumerate}

Positive solutions are sought in the product cone $K = K_1 \times K_2 \subset C(I) \times C(I)$, where for each $i \in \{1, 2\}$:
\begin{equation}\label{eq:cone_Ki}
	K_i := \left\{ w \in P_i :  \min_{t \in [a, b]} w(t) \ge c_i \|w\|_\infty \right\},
\end{equation}
with $[a, b] \subseteq I$, $c_i \in (0, 1]$ and $\tilde{c}_i\in[0,1]$ given by $(H_2)$, and $P_i$ denotes the cone
\begin{equation*}\label{eq:cone_Pi}
	P_i := \left\{ w \in C(I) : \min_{t \in I} w(t) \ge \tilde{c}_i \|w\|_\infty \right\}.
\end{equation*}
The inequalities within cones of the form \eqref{eq:cone_Ki} are known as Harnack inequalities, while $c_i$ and $\tilde c_i$ are the corresponding Harnack constants. It is important to emphasize the qualitative distinction between these two constants: $c_i$ is strictly positive, whereas $\tilde{c}_i$ is only nonnegative and may vanish, i.e., $\tilde{c}_i=0$. The latter situation arises naturally in boundary value problems with homogeneous Dirichlet conditions at the endpoints of $[A,B]$, since admissible solutions may vanish at $t=A$ or $t=B$, yielding $\min_{t\in[A,B]}w(t)=0$. Hence, a positive global Harnack constant cannot generally be expected in this setting.

In contrast, for periodic problems, as well as for problems with strictly positive Neumann or Robin boundary conditions, the associated Green's functions are typically strictly positive on $I\times I$, allowing one to take $\tilde{c}_i>0$. Nevertheless, $\tilde{c}_i=0$ may still occur in periodic settings if the Green's function vanishes somewhere in $I\times I$, as may happen for certain higher-order periodic operators or when the operator lies at the boundary of the regime in which a maximum principle holds \cite{CabadaBook}.

Cones of this type have been previously considered in \cite{web}.

For each $i, j \in \{1, 2\}$, we define the local minima and global maxima of the kernels by
\begin{equation*}\label{eq:kernel_extrema}
 m_j^i := \min_{(t, s) \in [a, b] \times [a, b]} G_j^i(t, s), \qquad M_j^i := \max_{(t, s) \in I \times I} G_j^i(t, s).
\end{equation*}
 Note that condition $(H_2)$--$(iii)$ implies that $m_i^i>0$ for $i\in\{1,2\}$.
Moreover, we introduce the notation
\begin{align*}
	\underline{\lambda} &:= \min\{\lambda_1, \lambda_2\} > 0, \qquad \bar{\lambda} := \max\{\lambda_1, \lambda_2\} > 0,  \\[1ex]
	m^i &:= \min\{m_1^i, m_2^i\} , \qquad M^i := \max\{M_1^i, M_2^i\} .
\end{align*}

To prove our existence result, we apply the vectorial version of the
Krasnosel'ski\u{\i}--Guo theorem (Theorem~\ref{th1.1}) to the operator $\mathcal{T}:K:=K_1\times K_2\to K_1\times K_2$ given by $\mathcal{T}=(\mathcal{T}_1,\mathcal{T}_2)$, with $\mathcal{T}_i:K\to K_i$, defined by
\[\mathcal{T}_i(u,v)(t):= \lambda_1\int_{I} G_1^i(t,s)f_1(s,u(s),v(s))\,ds
+ \lambda_2\int_{I} G_2^i(t,s)f_2(s,u(s),v(s))\,ds,\quad t\in I, \quad  i=1,2.\]

Under assumptions $(H_1)$--$(H_3)$, we have $\mathcal{T}(K)\subseteq K$. Indeed, let $(u,v)\in K$ and set
\begin{equation}
	\label{2.14}
	F_j(s):=\lambda_j f_j(s,u(s),v(s))\geq0,
	\qquad s\in I,\quad j=1,2.
\end{equation}
Since $G_j^i(t,s)\geq0$, we have $\mathcal{T}_i(u,v)(t)\geq0$ for all $t\in I$.
Moreover, by $(H_1)$ and $(H_3)$, $\mathcal{T}_i(u,v)\in C(I)$. Applying the upper bound \eqref{eq:H2_upper} from $(H_2)$, for every $t \in I$
\[
\mathcal{T}_i(u, v)(t) = \sum_{j=1}^2 \int_{I} G_j^i(t, s) F_j(s) \, ds \le \sum_{j=1}^2 \int_{I} \Phi_j^i(s) F_j(s) \, ds,
\]
hence it follows that
\begin{equation}\label{eq:norm_bound}
	\|\mathcal{T}_i(u, v)\|_\infty = \max_{t \in I} \mathcal{T}_i(u, v)(t) \le \sum_{j=1}^2 \int_{I} \Phi_j^i(s) F_j(s) \, ds.
\end{equation}
Note that, here, if $G_j^i\equiv0$ for $j\neq i$, then $\Phi_j^i\equiv0$, so the corresponding
term vanishes automatically.

On the other hand, \eqref{eq:H2_lower_loc} and \eqref{eq:H2_lower_glob} from $(H_2)$ imply
\[
\mathcal{T}_i(u,v)(t)
\geq
c_i\sum_{j=1}^2\int_{I}\Phi_j^i(s)F_j(s)\,ds
\geq c_i\|T_i(u,v)\|_\infty,
\qquad t\in[a,b],
\]
and
\[
\mathcal{T}_i(u,v)(t)
\geq
\tilde c_i\sum_{j=1}^2\int_{I}\Phi_j^i(s)F_j(s)\,ds
\geq \tilde c_i\|T_i(u,v)\|_\infty,
\qquad t\in I.
\]
Therefore, $\mathcal{T}_i(u,v)\in K_i$ for $i=1,2$.

Moreover, $\mathcal{T}$ is completely continuous, that is, continuous and maps bounded sets into relatively compact ones, as a direct
consequence of the Ascoli--Arzel\`a theorem.

Therefore, by Theorem~\ref{th1.1}, obtaining a solution with prescribed bounds on the norms of its components reduces to verifying, for each $i\in\{1,2\}$, the compression or expansion inequalities for $\mathcal{T}_i$ on $\left\{x\in \overline{K}_{r,R}: \|x_i\| = r_i\right\}$ and $\left\{x\in \overline{K}_{r,R}: \|x_i\| = R_i\right\}$. This is done in Sections~\ref{sec:main_results} and~\ref{sec:asymptotic_multiplicity}.

To conclude this section, we recall the fundamental properties of the fixed point index that will be used throughout this work. For further details, we refer the reader to \cite{amann, guolak}.

\begin{proposition}\label{prop_index}
	Let $K$ be a cone in a normed space, let $U\subset K$ be a bounded relatively open set, and let $\mathcal{T}:\overline{U}\to K$ be a compact map with no fixed points on the boundary of $U$. Then the fixed point index of $\mathcal{T}$ on $U$ with respect to $K$, denoted by $i_K(\mathcal{T},U)$, satisfies the following properties:
	\begin{enumerate}[$(i)$]
		\item \textup{(Additivity)} Suppose that $U_1$ and $U_2$ are disjoint relatively open subsets of $U$. If $0\notin (I-\mathcal{T})\bigl(\overline{U}\setminus(U_1\cup U_2)\bigr)$ then
		\[i_K(\mathcal{T},U)=i_K(\mathcal{T},U_1)+i_K(\mathcal{T},U_2).\]
		
		\item \textup{(Existence)} If $i_K(\mathcal{T},U)\neq 0$, then there exists $x\in U$ such that $\mathcal{T}x=x$.
	\end{enumerate}
\end{proposition}

The following fixed point index computation follows from \cite[Proposition~2.2]{LFP_JRL} and will be useful in Section~\ref{sec:multiplicity}.

\begin{proposition}
	\label{prop2.2}
	Let $R\in\mathbb{R}^2$ with $R_i>0$ for $i=1,2$, and define
	\[K_R:=\left\{x=(x_1,x_2)\in K:\|x_i\|<R_i,\ i=1,2\right\};\]
	\[\overline{K}_R:=\left\{x=(x_1,x_2)\in K:\|x_i\|\leq R_i,\ i=1,2\right\}.\]
	Assume that $\mathcal{T}=(\mathcal{T}_1,\mathcal{T}_2):
	\overline{K}_R\to K$ is a compact map satisfying $\|\mathcal{T}_i x\|> \|x_i\|$ whenever $\|x_i\|=R_i$, $i=1,2$.
	Then 
	\[i_K(\mathcal{T},K_R)=0.\]
\end{proposition}

\section{Existence and non-existence results }
\label{sec:main_results}

In this section, we establish our main existence and non-existence results for solutions with nontrivial components to problem~\eqref{1.1}.

We first state our main existence result.

\begin{theorem}
	\label{th2.1}
	Assume that conditions $(H_1)\text{--}(H_3)$ hold and, moreover, for each $i\in\{1,2\}$ there exist $\alpha_i,\beta_i>0$ with $\alpha_i\neq\beta_i$ such that the following
	conditions are fulfilled:
	\begin{equation}
		\label{2.5}
		\lambda_1m_1^i f_1(t,x_1,x_2)+\lambda_2m_2^i f_2(t,x_1,x_2)\geq \frac{\alpha_i}{b-a},
		\, \text{for } x_i\in[c_i\alpha_i,\alpha_i],\ 
		x_j\in[c_jr_j,R_j]\ (j\neq i),\ t\in[a,b]; \text{ and }
	\end{equation}
	\begin{equation}
		\label{2.6}
		\lambda_1M_1^i f_1(t,x_1,x_2)+\lambda_2M_2^i f_2(t,x_1,x_2)\leq \frac{\beta_i}{B-A},
		\,\text{for } x_i\in[\tilde{c}_i\beta_i,\beta_i],\ 
		x_j\in[\tilde{c}_jr_j,R_j]\ (j\neq i),\ t\in I,
	\end{equation}
	where $	r_i:=\min\{\alpha_i,\beta_i\}$ and $	R_i:=\max\{\alpha_i,\beta_i\}$, $i=1,2$.
	
	Then problem \eqref{1.1} has at least one positive solution $(u,v)$ with $r_1\leq \|u\|_\infty\leq R_1$ and $	r_2\leq \|v\|_\infty\leq R_2$.
\end{theorem}

\noindent
{\bf Proof.}
	Let us verify that $\mathcal{T}$ satisfies the assumptions of
	Theorem~\ref{th1.1} and, consequently, has a fixed point
	$(x,y)\in \overline{K}_{r,R}$, which is a solution of system~\eqref{1.1}.
	
	First, for each $i\in\{1,2\}$, we prove that $\|\mathcal{T}_i(x_1,x_2)\|_\infty\geq\|x_i\|_\infty$ whenever $\|x_i\|_\infty=\alpha_i.$ Indeed, if $(x_1,x_2)\in \overline{K}_{r,R}$ and
	$\|x_i\|_\infty=\alpha_i$, then $c_i\alpha_i\leq x_i(t)\leq\alpha_i$ and $	c_jr_j\leq x_j(t)\leq R_j,\, j\neq i,$ for all $t\in[a,b]$. Hence, by \eqref{2.5}, for each $t\in[a,b]$, we have
	\begin{equation}
		\label{2.7}
	\begin{aligned}
		\mathcal{T}_i(x_1,x_2)(t)=\sum_{j=1}^2 \int_{I} G_j^i(t, s) F_j(s) \, ds &\ge \sum_{j=1}^2 \int_a^b G_j^i(t, s) F_j(s) \, ds\\
		& \ge \sum_{j=1}^2 \int_a^b m_j^i F_j(s) \, ds \geq \frac{\alpha_i}{b-a}\int_a^b \,ds=\alpha_i.
	\end{aligned}
	\end{equation}
Thus, taking the supremum over $t\in [a, b]$, we obtain $\|\mathcal{T}_i(x_1,x_2)\|_\infty\geq\|x_i\|_\infty$, as desired.
	
	Second, for each $i\in\{1,2\}$, we show that $	\|\mathcal{T}_i(x_1,x_2)\|_\infty\leq\|x_i\|_\infty$
	whenever $\|x_i\|_\infty=\beta_i$. Indeed, if
	$(x_1,x_2)\in \overline{K}_{r,R}$ and $\|x_i\|_\infty=\beta_i$, then $	\tilde{c}_i\beta_i\leq x_i(t)\leq\beta_i$ and $\tilde{c}_jr_j\leq x_j(t)\leq R_j,
	\, j\neq i$ for all $t\in I$. Thus, by \eqref{2.6}, for every $t\in I$,
\[
\begin{aligned}
	\mathcal{T}_i(x_1,x_2)(t)=\sum_{j=1}^2 \int_{I} G_j^i(t, s) F_j(s) \, ds \le \sum_{j=1}^2 \int_{I} M_j^i F_j(s) \, ds \le \frac{\beta_i}{B-A}\int_{I} \,ds=\beta_i.
\end{aligned}
\]
Hence, $\|\mathcal{T}_i(x_1,x_2)\|_\infty\leq\|x_i\|_\infty.$

The conclusion now follows from Theorem~\ref{th1.1}. Indeed,
$\mathcal{T}_i$ is compressive when $\alpha_i<\beta_i$ and expansive
when $\alpha_i>\beta_i$. \qed

\begin{remark}
Note that, when $m_j^i>0$ for $j=1,2$, it is sufficient to assume
\begin{equation*}
	\begin{aligned}
		\lambda_1m_1^i f_1(t,x_1,x_2)\geq \frac{ \alpha_i}{b-a},
		\quad &\text{for } x_i\in[c_i\alpha_i,\alpha_i],\ 
		x_j\in[c_jr_j,R_j]\ (j\neq i),\ t\in[a,b]; \quad\text{or}\\
		\lambda_2m_2^i f_2(t,x_1,x_2)\geq \frac{\alpha_i}{b-a},
		\quad &\text{for } x_i\in[c_i\alpha_i,\alpha_i],\ 
		x_j\in[c_jr_j,R_j]\ (j\neq i),\ t\in[a,b],
	\end{aligned}
\end{equation*}
instead of condition \eqref{2.5} in Theorem~\ref{th2.1}.
Indeed, this follows from the nonnegativity of the nonlinearities and kernels,
which ensures that the last inequality in \eqref{2.7} remains valid.
\end{remark}

\begin{remark}\label{rmak_tildec}
		The role of the constants $\tilde{c}_i$ in hypothesis \eqref{2.6} deserves special attention. If $\tilde{c}_i>0$, the upper bound condition is imposed only on a truncated box whose $i$-th component is bounded away from zero, namely on the interval $[\tilde{c}_i\beta_i,\beta_i]$. Therefore, no control of the nonlinearities is needed near $x_i=0$. Consequently, the theorem applies to a wider class of systems, including nonlinearities that may be singular at the origin. On the other hand, if $\tilde{c}_i=0$, then the estimate must hold on the whole interval $[0,\beta_i]$, which is a stronger requirement.
		
		The above mentioned difference between the cases $\tilde{c}_i>0$ and $\tilde{c}_i=0$ motivates the study of distinct asymptotic conditions in Section~\ref{sec:asymptotic_multiplicity} depending on which case applies.  
\end{remark}

\begin{remark}
Theorem~\ref{th2.1} extends \cite[Theorem~4.2]{LFP_JRL}. Indeed, when $[a,b]=I$, $\lambda_1=\lambda_2=1$, and $G_2^1(t,s)=G_1^2(t,s)=0$ for all $(t,s)\in I\times I$, problem~\eqref{1.1} reduces to \eqref{1.4}, and conditions \eqref{2.5} and \eqref{2.6} in Theorem~\ref{th2.1} improve assumptions (13) and (14) in \cite[Theorem~4.2]{LFP_JRL}. Note that \cite[Theorem~4.2]{LFP_JRL} concerns systems of second-order equations with periodic boundary conditions, but it can be easily adapted to systems of Hammerstein equations of the form \eqref{1.4} with strictly positive kernels.
\end{remark}

The following result provides conditions under which the components of a solution are trivial, that is, it proves non-existence of positive coexistence solutions.

\begin{theorem}\label{th:non_existence}
	Assume that conditions $(H_1)\text{--}(H_3)$ hold. Moreover, fix $i \in \{1, 2\}$ and assume that one of the following conditions is fulfilled:
	\begin{enumerate}[$(a)$]
		\item For all $t \in I$, $x_i > 0$, and $x_j \ge 0$ with $j \neq i$,
		\begin{equation}\label{eq:cond_sub_global}
			\lambda_1 M_1^i f_1(t, x_1, x_2) + \lambda_2 M_2^i f_2(t, x_1, x_2) < \frac{x_i}{B-A}; 
		\end{equation}
		\item For all $t \in [a,b]$, $x_i > 0$, and $x_j \ge 0$ with $j \neq i$,
		\begin{equation}\label{eq:cond_super_global}
			\lambda_1 m_1^i f_1(t, x_1, x_2) + \lambda_2 m_2^i f_2(t, x_1, x_2) > \frac{x_i}{b - a}.
		\end{equation}
	\end{enumerate}
	Then, if $(u_1,u_2)\in K$ is a solution for problem \eqref{1.1}, then $u_i(t)=0$ for all $t\in I$.
\end{theorem}

\noindent \textit{Proof.} Let $(u_1, u_2) \in K$ be a solution of \eqref{1.1}. First, suppose that condition \eqref{eq:cond_sub_global} holds. Arguing by contradiction, assume that $\|u_i\|_\infty > 0$, and let $t_0 \in I$ be such that $u_i(t_0) = \|u_i\|_\infty$. 
	
	Since $u_i \in K_i$, we have $u_i(t) \ge 0$ for all $t\in I$ and, moreover, $u_i(t)>0$ for all $t\in [a, b]$. Thus, in view of \eqref{eq:cond_sub_global}, we obtain
	\begin{align*}
		\|u_i\|_\infty = u_i(t_0) 
		&\le \int_{I} \left[ \lambda_1 M_1^i f_1\big(s, u_1(s), u_2(s)\big) + \lambda_2 M_2^i f_2\big(s, u_1(s), u_2(s)\big) \right] ds \\
		&< \int_{I} \frac{u_i(s)}{B-A} \, ds \le \int_{I} \frac{\|u_i\|_\infty}{B-A} \, ds = \|u_i\|_\infty,
	\end{align*}
	which is a contradiction. Consequently, $u_i(t)=0$ for all $t\in I$.
	
	Next, suppose that condition \eqref{eq:cond_super_global} holds. Again proceeding by contradiction, assume that $\|u_i\|_\infty > 0$. Let $\eta_i \in [a, b]$ satisfy
	\[
	u_i(\eta_i) = \min_{t \in [a, b]} u_i(t) \ge c_i \|u_i\|_\infty > 0.
	\]
	In particular, $u_i(t) \ge u_i(\eta_i) > 0$ for all $t \in [a, b]$. Applying \eqref{eq:cond_super_global}, for every $t \in [a, b]$ we obtain
	\begin{align*}
		u_i(t) &\ge \int_a^b \left[ \lambda_1 m_1^i f_1\big(s, u_1(s), u_2(s)\big) + \lambda_2 m_2^i f_2\big(s, u_1(s), u_2(s)\big) \right] ds \\
		&> \int_a^b \frac{u_i(s)}{b - a} \, ds \ge \int_a^b \frac{u_i(\eta_i)}{b - a} \, ds = u_i(\eta_i).
	\end{align*}
	Taking the minimum over $t \in [a, b]$ on the left-hand side yields $u_i(\eta_i) > u_i(\eta_i)$, which is an absurdity. Hence, $u_i(t)=0$ for all $t\in I$. \qed
	
The preceding result leads directly to the following non-existence theorem for nontrivial solutions.
	
\begin{corollary}\label{cor:non_existence}
	Assume that conditions $(H_1)\text{--}(H_3)$ hold.
	\begin{enumerate}[$(i)$]
		\item If, for each $i\in\{1,2\}$, one of conditions $(a)$
		or $(b)$ in Theorem~\ref{th:non_existence} is satisfied, then
		problem~\eqref{1.1} admits no nontrivial solution in $K$.
		
		\item Suppose that, for some $i\in\{1,2\}$, one of
		conditions $(a)$ or $(b)$ in Theorem~\ref{th:non_existence} holds.
		Let $j\in\{1,2\}$, $j\neq i$, and assume that
		\[
		f_1(t,x_1,x_2)=f_2(t,x_1,x_2)=0
		\qquad\text{whenever }x_i=0,\quad
		\text{for all }x_j\geq0\text{ and }t\in I.
		\]
		Then problem~\eqref{1.1} admits no nontrivial solution in $K$.
	\end{enumerate}
\end{corollary}

\noindent \textit{Proof.} Statement $(i)$ follows directly from the preceding theorem. 

For proving $(ii)$, suppose that $(u_1,u_2)$ is a solution for problem \eqref{1.1} and that for $i=1$ condition $(a)$ or $(b)$ in Theorem~\ref{th:non_existence} holds, which implies that $u_1(t)=0$ for all $t\in I$. Then
\[u_2(t)=\lambda_1\int_{I} G_1^2(t,s)f_1(s,0,u_2(s))\,ds
+ \lambda_2 \int_{I} G_2^2(t,s)f_2(s,0,u_2(s))\,ds, \quad t\in I, \]
and, since $f_1(t,0,x_2)=0=f_2(t,0,x_2)$ for all $x_2\geq 0$ and $t\in I$, it follows that $u_2(t)=0$ for all $t\in I$.  
\qed

In what follows, we illustrate the applicability of the previous results by establishing existence and non-existence of radially symmetric solutions to an elliptic system on an annulus subject to coupled boundary conditions.
	
		\begin{example}
			\label{ej:2.1}
Consider the elliptic system
\begin{equation}\label{eq:pde_system}
	\left\{
	\begin{aligned}
		-\Delta U &=
		\frac{4\lambda}{|x|^4}
		\tanh(U)\bigl(1+\tanh(V)\bigr),
		&& 1< |x|<2,\\[0.6em]
		-\Delta V &=
		\frac{4\lambda}{|x|^4}
		\tanh(U)\tanh(V),
		&& 1< |x|<2,
	\end{aligned}
	\right.
\end{equation}
in the annular domain
$\Omega:=\{x\in\mathbb{R}^3:1<|x|<2\}$, subject to the coupled
Robin boundary conditions
\begin{equation}\label{eq:pde_bc1} \left\{ \begin{aligned} U(x) - \frac{1}{2} \frac{\partial U}{\partial \nu}(x) - V(x) &= 0, \quad |x| = 1; \qquad & U(x) + 2\frac{\partial U}{\partial \nu}(x) &= 0, \quad |x| = 2, \\[0.8em] V(x) - \frac{1}{2} \frac{\partial V}{\partial \nu}(x) - U(x) &= 0, \quad |x| = 1; \qquad & V(x) + 2 \frac{\partial V}{\partial \nu}(x) &= 0, \quad |x| = 2. \end{aligned} \right. \end{equation}

We show that this problem admits at least one radially symmetric
solution $(U,V)$---that is, depending only on $r=|x|$--- with both components nontrivial whenever
$\lambda>3.0975$. In contrast, we claim that, for
$\lambda<6/17\approx0.3529$, the trivial solution is the only
radially symmetric solution.

	In order to establish the existence of radial solutions $U(|x|)=u(r)$ and $V(|x|)=v(r)$, system \eqref{eq:pde_system} can be rewritten as
			\begin{equation*}
				\left\{
				\begin{aligned}
					-\dfrac{1}{r^{2}}(r^2 u')' &=
					\frac{4\lambda}{r^4}
					\tanh(u)\bigl(1+\tanh(v)\bigr),
					&& 1< r<2,\\[0.6em]
					-\dfrac{1}{r^{2}}(r^2 v')' &=
					\frac{4\lambda}{r^4}
					\tanh(u)\tanh(v),
					&& 1< r<2,
				\end{aligned}
				\right.
			\end{equation*}
			and the boundary conditions \eqref{eq:pde_bc1} as
			\begin{equation} \left\{ \begin{aligned} u(1) - \frac{1}{2} u'(1) - v(1) &= 0,  \qquad & u(2) + 2u'(2) &= 0,  \\[0.8em] v(1) - \frac{1}{2} v'(1) - u(1) &= 0, \qquad & v(2) + 2 v'(2) &= 0. \end{aligned} \right. \end{equation}
	Set now $u(t)=u(t(r))$ and $v(t)=v(t(r))$, with $t(r):=3-2/r$ for $r\in [1,2]$. 		
	Then it is straightforward to verify that the radially symmetric solutions of \eqref{eq:pde_system}--\eqref{eq:pde_bc1} correspond to the solutions of the system of ordinary differential equations
		\begin{equation}\label{eq:ode_system}
			\left\{
			\begin{aligned}
				-u''(t) &= \lambda f_1\big(u(t), v(t)\big), \quad t \in I := [1, 2], \\[0.6em]
				-v''(t) &= \lambda f_2\big(u(t), v(t)\big), \quad t \in I, \\[0.6em]
				u(1) + u'(1) - & v(1) = 0, \qquad u(2) + u'(2) = 0, \\[0.6em]
				v(1) + v'(1) - & u(1) = 0, \qquad v(2) + v'(2) = 0.
			\end{aligned}
			\right.
		\end{equation}
		with the nonlinearities $f_1,f_2:[0,\infty)^2\to [0,\infty)$ given by 
		\[f_1(x_1,x_2):=\tanh{(x_1)}(1+\tanh{(x_2)}) \quad \text{ and } \quad f_2(x_1,x_2):=\tanh{(x_1)} \tanh{(x_2)}.\]
	This problem is a particular case of \eqref{1.2}. Moreover, following \cite{Mou}, one can establish the equivalence between the solutions of \eqref{eq:ode_system} and those of the integral system \eqref{1.1}, in which the four Green's functions on $I^2$ are given by
	\begin{equation*}
		G_1^1(r,s)=G_2^2(r,s)
		=\frac{\bigl(3-\max\{r,s\}\bigr)\min\{r,s\}}{3},
		\qquad
		G_2^1(r,s)=G_1^2(r,s)
		=\frac{(3-r)(3-s)}{3}.
	\end{equation*}
Taking $[a,b]=I$, the Green's functions $G_j^i$ are continuous and strictly positive on $I^2$. Moreover, they satisfy conditions $(H_1)$ and $(H_2)$ with $ m_j^i=1/3,\, M_2^1=M_1^2=4/3,\, M_1^1=M_2^2=3/4,\, i,j\in\{1,2\},$ while the corresponding cone constants are $c_1=c_2=\tilde c_1=\tilde c_2=1/2>0$. See Appendices~A, B, and~C for the details.
	
	Before proving the claim, notice that the nonlinearities $f_1, f_2$ are strictly increasing with respect to each variable and satisfy
	\begin{equation}
		\label{2.25}
	f_1(x_1, x_2) \le 2 x_1, \qquad f_2(x_1, x_2) \le x_1, \qquad \text{for all } x_1, x_2 \ge 0,
	\end{equation}
	as well as
	\begin{equation}
		\label{2.26}
	f_1(x_1, x_2) \le 2, \qquad f_2(x_1, x_2) \le 1, \qquad \text{for all } x_1, x_2 \ge 0.
	\end{equation}
	Moreover, for any fixed $\alpha_1, \alpha_2 > 0$, we have
\begin{equation*}
	\label{2.27}
(f_1 + f_2)(x_1, x_2) \ge \tanh{\left(\frac{\alpha_1}{2}\right)}\left[1 + 2\tanh{\left(\frac{\alpha_2}{2}\right)}\right] \quad \text{for all } x_1 \ge \frac{\alpha_1}{2} \text{ and } x_2 \ge \frac{\alpha_2}{2}.
\end{equation*}

To prove non-existence, let $\lambda<6/17$. For $x_1>0$ and $x_2\geq0$, estimate~\eqref{2.25} gives
\[
\lambda M_1^1 f_1(x_1, x_2) + \lambda M_2^1 f_2(x_1, x_2) \le \lambda \left[ \frac{3}{4} (2x_1) + \frac{4}{3} x_1 \right] = \frac{17}{6} \lambda x_1.
\]
Thus, condition~(a) in Theorem~\ref{th:non_existence} holds, and hence $x_1(t)=0$ for all $t\in I$ and every solution $(x_1,x_2)\in K$. Moreover, $f_1(0,x_2)=f_2(0,x_2)=0$ for all $x_2\geq0$, so Corollary~\ref{cor:non_existence}(ii) yields $x_2(t)=0$ for all $t\in I$. Therefore, problem~\eqref{eq:ode_system} admits no nontrivial solution in $K$ for $\lambda<6/17$.

We now verify conditions \eqref{2.5} and \eqref{2.6} of Theorem~\ref{th2.1}, choosing $\alpha_i$ and $\beta_i$ so as to minimize the resulting lower bound for $\lambda$. First, observe that
\begin{align*}
	\lambda M_1^1 f_1(x_1, x_2) + \lambda M_2^1 f_2(x_1, x_2) &\le \lambda\left(\frac{3}{4}\cdot 2 + \frac{4}{3}\cdot 1\right) = \frac{17}{6}\lambda \le \beta_1, \quad \text{for } (x_1, x_2) \in [\beta_1/2, \beta_1] \times [\alpha_2/2, \beta_2], \\[0.5em]
	\lambda M_1^2 f_1(x_1, x_2) + \lambda M_2^2 f_2(x_1, x_2) &\le \lambda\left(\frac{4}{3}\cdot 2 + \frac{3}{4}\cdot 1\right) = \frac{41}{12}\lambda \le \beta_2, \quad \text{for } (x_1, x_2) \in [\alpha_1/2, \beta_1] \times [\beta_2/2, \beta_2].
\end{align*}
Accordingly, we set $\beta_1 := \frac{17}{6}\lambda$ and $\beta_2 := \frac{41}{12}\lambda$.

On the other hand, condition \eqref{2.5} requires that for every $(x_1, x_2) \in \left([\alpha_1/2, \alpha_1] \times [\alpha_2/2, \beta_2]\right) \cup ([\alpha_1/2, \beta_1] \times [\alpha_2/2, \alpha_2])$,
\[
\lambda\bigl(m_1^i f_1 + m_2^i f_2\bigr)(x_1, x_2) = \frac{\lambda}{3}\bigl(f_1 + f_2\bigr)(x_1, x_2) \ge \frac{\lambda}{3}\tanh{\left(\frac{\alpha_1}{2}\right)}\left[1 + 2\tanh{\left(\frac{\alpha_2}{2}\right)}\right] \ge \alpha_i,
\]
which leads to the requirement
\[
\lambda \ge \frac{3\alpha_i}{\tanh(\alpha_1/2)\left[1 + 2\tanh(\alpha_2/2)\right]}.
\]
Taking $\alpha_1 = \alpha_2 = \alpha > 0$, both inequalities reduce to
\[
\lambda \ge \Phi_{\mathrm{const}}(\alpha) := \frac{3\alpha}{\tanh(\alpha/2)\left[1 + 2\tanh(\alpha/2)\right]}.
\]
The function $\Phi_{\mathrm{const}}$ attains its unique global minimum at $\alpha^* \approx 1.7357$, with $\Phi_{\mathrm{const}}(1.7357) \approx 3.0975$. Thus, for any $\lambda \ge 3.0975$, we obtain
\[
\beta_1 = \frac{17}{6}\lambda \ge 8.7762 > 1.7357 = \alpha_1, \qquad \beta_2 = \frac{41}{12}\lambda \ge 10.5831 > 1.7357 = \alpha_2.
\]
Therefore, Theorem~\ref{th2.1} yields at least one positive solution $(u,v)$ of problem~\eqref{eq:ode_system} satisfying $\alpha_1 \le \|u\|_\infty \le \beta_1$ and $\alpha_2 \le \|v\|_\infty \le \beta_2$.
		\end{example}

\section{On asymptotic criteria}
\label{sec:asymptotic_multiplicity}

While the application of Theorem~\ref{th2.1} requires explicit bounds for the Green's functions, existence results (e.g., \cite{doOLU,doOLU2,wang,Y,JRL}) often rely only on their positivity together with asymptotic conditions on the nonlinearities that are more readily interpretable than those in the theorem. Motivated by this observation, and building on this result, we develop new existence results for positive solutions with nontrivial components to problem~\eqref{1.1}, based on directional and diagonal asymptotic growth conditions. 

We begin by introducing some notation that will be used throughout the section.
For $i\in\{1,2\}$, given a continuous function $g:I\times[0,\infty)^2\to [0,\infty)$, we define
\[(g)_0^i:=\lim_{x_i\to0^+}\frac{g(t,x_1,x_2)}{x_i},\qquad(g)_\infty^i:=\lim_{x_i\to+\infty}\frac{g(t,x_1,x_2)}{x_i}.\]

Moreover, for each $i \in \{1,2\}$, let
\[
\eta_j^i := \begin{cases} 1, & m_j^i > 0, \\ 0, & m_j^i = 0, \end{cases}
\qquad \text{and} \qquad
\delta_j^i := \begin{cases} 1, & M_j^i > 0, \\ 0, & M_j^i = 0, \end{cases}
\quad (j=1,2),
\]
and set 
\[S_i^\eta := \eta_1^i f_1 + \eta_2^i f_2 \quad \text{ and } \quad S_i^\delta := \delta_1^i f_1 + \delta_2^i f_2.\] 
Notice that condition $(H_2)\text{--}(iii)$ implies that $m_i^i>0$ for $i\in\{1,2\}$ and so $\eta_i^i=\delta_i^i=1$ for $i\in\{1,2\}$. Hence, the functions $S_i^\eta$ and $S_i^\delta$ are, in fact, either equal to $f_i$ or to $f_1+f_2$, depending on the values of $\eta_j^i$ and $\delta_j^i$, $j\neq i$, respectively. 
Our existence results below rely on the asymptotic behavior of these functions $S_i^\eta$ and $S_i^\delta$ at zero and at infinity.

For a fixed pair of distinct indices $(k,l) \in \{(1,2), (2,1)\}$, we consider the following conditions:
\begin{enumerate}
	\item[($\mathcal{A}_k$)] $(S_k^\eta)_0^k = +\infty$ and $(S_k^\delta)_\infty^k = 0$, uniformly for $(t,x_l) \in I \times \mathbb{R}^+$;
	\item[($\mathcal{B}_k$)] $(S_k^\delta)_0^k = 0$ and $(S_k^\eta)_\infty^k = +\infty$, uniformly for $(t,x_l) \in I \times \mathbb{R}^+$;
	\item[($\widetilde{\mathcal{A}}^0_k$)] $(S_k^\eta)_0^k = +\infty$ uniformly for $(t,x_l) \in I \times J$ for each compact $J \subset \mathbb{R}_+:=(0,\infty)$ and $(S_k^\delta)_\infty^k = 0$, uniformly for $(t,x_l) \in I \times J$ for each compact $J \subset \mathbb{R}^+$;
	\item[($\widetilde{\mathcal{B}}^0_k$)] $(S_k^\delta)_0^k = 0$ uniformly for $(t,x_l) \in I \times J$ for each compact $J \subset \mathbb{R}^+$, and $(S_k^\eta)_\infty^k = +\infty$, uniformly for $(t,x_l) \in I \times J$ for each compact $J \subset \mathbb{R}_+$.
	\item[($\mathcal{C}_k$)] for each $R > 0$,
	\[
	M_k^*(R) := \sup_{t\in I, x_k\in[0,R], x_l\in [0,\infty)} S_k^\delta(t,x_1,x_2) < \infty;
	\]
\end{enumerate}

Under conditions $(\mathcal{A}_k)$ and $(\widetilde{\mathcal{A}}^0_k)$, we say that $S_k^\eta$ exhibits \textit{superlinear behavior at zero} and $S_k^\delta$ exhibits \textit{sublinear behavior at infinity} in the $x_k$-direction, uniformly with respect to $t$ and $x_l$. Conversely, under conditions $(\mathcal{B}_k)$ and $(\widetilde{\mathcal{B}}^0_k)$, we say that $S_k^\delta$ exhibits \textit{sublinear behavior at zero} and $S_k^\eta$ exhibits \textit{superlinear behavior at infinity} in the $x_k$-direction, uniformly with respect to $t$ and $x_l$.

We now state the first existence result for solutions with both components nontrivial under uniform asymptotic conditions of this type. Let us focus on the case in which the constant $\tilde{c}_k$ in assumption $(H_2)$ is zero.

\begin{theorem}
	\label{coro2.1}
	Assume that conditions $(H_1)\text{--}(H_3)$ hold. Let $(k,l)\in\{(1,2),(2,1)\}$ and suppose that $\tilde{c}_k=0$, $\tilde{c}_l\geq0$, and either:
		\begin{enumerate}
			\item[\textup{$(a)$}] the $k$-th component satisfies $(\mathcal{A}_k)$ and $(\mathcal{C}_k)$, while the $l$-th component satisfies $(\widetilde{\mathcal{A}}_l^0)$;
			\item[\textup{$(b)$}] $M_k^l=M_l^k=0$, and either the $k$-th component satisfies $(\mathcal{B}_k)$, while the $l$-th component satisfies $(\widetilde{\mathcal{A}}_l^0)$ or $(\widetilde{\mathcal{B}}_l^0)$; or  the $k$-th component satisfies $(\mathcal{A}_k)$ and $(\mathcal{C}_k)$, while the $l$-th component satisfies $(\widetilde{\mathcal{A}}_l^0)$ or $(\widetilde{\mathcal{B}}_l^0)$.
		\end{enumerate}
	Then, for every $\lambda_1,\lambda_2>0$, problem~\eqref{1.1} admits at least one positive solution $(u,v)$ such that $\|u\|_\infty>0$ and $\|v\|_\infty>0$.
\end{theorem}

\noindent
{\bf Proof.} Let us prove $(a)$. Fix $\lambda_1,\lambda_2>0$. Assume that $\tilde c_1=\tilde c_2=0$ and that $(\mathcal{A}_1)$, $(\mathcal{C}_1)$, and $(\widetilde{\mathcal{A}}_2^0)$ hold, with $m_j^i,M_j^i>0$ for $i,j\in\{1,2\}$. The remaining configurations can be treated analogously.

Since $(f_1+f_2)_0^1=+\infty$ uniformly with respect to $t\in I$ and $x_2\in\mathbb{R}^+$, there exists $\alpha_1\equiv\alpha_1(\underline{\lambda})>0$ such that
\[(\lambda_1m_1^1f_1+\lambda_2m_2^1f_2)(t,x_1,x_2)\geq \underline{\lambda}\,m^1(f_1+f_2)(t,x_1,x_2)\geq \frac{x_1}{c_1(b-a)},\,\text{ for all } x_1\in(0,\alpha_1], \, x_2\in\mathbb{R}^+ \text{ and } t\in I. \]
Consequently,
\[
(\lambda_1m_1^1f_1+\lambda_2m_2^1f_2)(t,x_1,x_2)\geq\frac{c_1\alpha_1}{c_1(b-a)}=\frac{\alpha_1}{b-a}, \, \text{ for all } x_1\in[c_1\alpha_1,\alpha_1], \, x_2\in\mathbb{R}^+ \text{ and } t\in I.
\]
Thus, condition~\eqref{2.5} is satisfied for the first component.

On the other hand, the assumption $(f_1+f_2)_\infty^1=0$, uniformly with respect to $t\in I$ and $x_2\in\mathbb{R}^+$, guarantees the existence of $L_1\equiv L_1(\overline{\lambda})>\alpha_1$ sufficiently large such that
\[
\begin{aligned}
	(\lambda_1M_1^1f_1+\lambda_2M_2^1f_2)(t,x_1,x_2) &\leq \overline{\lambda}\,M^1(f_1+f_2)(t,x_1,x_2)\leq \frac{x_1}{B-A}, \\
	&\quad \text{for all } x_1\in[L_1,\infty), \, x_2\in\mathbb{R}_+ \text{ and } t\in I.
\end{aligned}
\]

Moreover, by $(\mathcal{C}_1)$,
\[
M_1^*(L_1)
:=
\sup_{(t,x_1,x_2)\in
	I\times[0,L_1]\times[0,\infty)}
(f_1+f_2)(t,x_1,x_2)
<\infty.
\]
Hence, choosing $\beta_1\equiv\beta_1(\overline{\lambda})\geq\max\left\{\overline{\lambda}M^1(B-A)M_1^*(L_1),L_1\right\}$, we obtain
\[(\lambda_1M_1^1f_1+\lambda_2M_2^1f_2)(t,x_1,x_2)\leq \overline{\lambda}\,M^1(f_1+f_2)(t,x_1,x_2)\leq \frac{\beta_1}{B-A}, \, \text{ for all } x_1\in[0,\beta_1], \, x_2\in\mathbb{R}_+ \text{ and } t\in I.\]
Therefore, inequality~\eqref{2.6} is also fulfilled for the first variable.

We now consider the second component. Since $(f_1+f_2)_0^2=+\infty$, uniformly with respect to $t\in I$ and $x_1\in J_1:=[c_1\alpha_1,\beta_1],$ there exists $\alpha_2\equiv\alpha_2(\underline{\lambda})>0$ such that
\[(\lambda_1m_1^2f_1+\lambda_2m_2^2f_2)(t,x_1,x_2)\geq \underline{\lambda}\,m^2(f_1+f_2)(t,x_1,x_2)\geq \frac{x_2}{c_2(b-a)}, \, \text{ for all } x_1\in J_1, \, x_2\in(0,\alpha_2] \text{ and } t\in I.
\]
In particular,
\[
(\lambda_1m_1^2f_1+\lambda_2m_2^2f_2)(t,x_1,x_2)
\geq
\frac{\alpha_2}{b-a}, \, \text{ for all } x_1\in J_1, \, x_2\in[c_2\alpha_2,\alpha_2] \text{ and } t\in I.
\]
Thus, condition~\eqref{2.5} is satisfied for the second variable.

Furthermore, since $(f_1+f_2)_\infty^2=0$ uniformly with respect to $t\in I$ and $x_1\in\widetilde J_1:=[0,\beta_1],$ there exists $L_2\equiv L_2(\overline{\lambda})>\alpha_2$ sufficiently large such that
\[(\lambda_1M_1^2f_1+\lambda_2M_2^2f_2)(t,x_1,x_2)\leq \overline{\lambda}\,M^2(f_1+f_2)(t,x_1,x_2)\leq \frac{x_2}{B-A}, \text{ for all } x_1\in\widetilde J_1, x_2\in[L_2,\infty) \text{ and } t\in I.
\]
Moreover, since the nonlinearities are continuous
\[
s(\beta_1,L_2):=\sup_{(t,x_1,x_2)\in
	I\times\widetilde J_1\times[0,L_2]} (f_1+f_2)(t,x_1,x_2) <\infty.
\]
It is therefore enough to choose $\beta_2\equiv\beta_2(\overline{\lambda})\geq \max\left\{\overline{\lambda}M^2(B-A)s(\beta_1,L_2),L_2\right\}$ so that
\[(\lambda_1M_1^2f_1+\lambda_2M_2^2f_2)(t,x_1,x_2)\leq \overline{\lambda}\,M^2(f_1+f_2)(t,x_1,x_2)\leq \frac{\beta_2}{B-A}, \text{ for all } x_1\in\widetilde J_1, x_2\in[0,\beta_2] \text{ and } t\in I.
\]
Hence, \eqref{2.6} holds for the second variable as well.

Note that, if $\tilde c_2>0$, the argument can be simplified. Indeed, we can directly choose a sufficiently large $\beta_2\equiv\beta_2(\overline{\lambda})>\alpha_2$, since $(f_1+f_2)_\infty^2=0$ uniformly with respect to $t\in I$ and $x_1\in\widetilde J_1$, such that
\[
\begin{aligned}
	(\lambda_1M_1^2f_1+\lambda_2M_2^2f_2)(t,x_1,x_2) &\leq \overline{\lambda}\,M^2(f_1+f_2)(t,x_1,x_2) \leq \frac{x_2}{B-A}, \\
	&\quad \text{for all } x_1\in\widetilde J_1, \, x_2\in[\tilde c_2\beta_2,\infty) \text{ and } t\in I,
\end{aligned}
\]
which implies \eqref{2.6} for $i=2$.

The remaining configurations in the decoupled case described in $(b)$ can be proven in a similar manner.
\qed

Let us consider now the case in which $\tilde{c}_k>0$. As already explained in Remark~\ref{rmak_tildec}, it implies that the estimate in \eqref{2.6} is less stringent, since no control of the nonlinearities is needed near $x_k=0$. Therefore, the asymptotic behavior described in $(\widetilde{\mathcal{A}}^0_k)$ and $(\widetilde{\mathcal{B}}^0_k)$ can be weakened imposing the assumptions uniformly only in compact intervals not including the zero. More precisely, for a fixed pair of distinct indices $(k,l) \in \{(1,2), (2,1)\}$, we consider the following conditions:
	\begin{enumerate}
			\item[($\widetilde{\mathcal{A}}_k$)] $(S_k^\eta)_0^k = +\infty$ and $(S_k^\delta)_\infty^k = 0$, uniformly for $(t,x_l) \in I \times J$ for each compact $J \subset (0,\infty)$;
			\item[($\widetilde{\mathcal{B}}_k$)] $(S_k^\delta)_0^k = 0$ and $(S_k^\eta)_\infty^k = +\infty$, uniformly for $(t,x_l) \in I \times J$ for each compact $J \subset (0,\infty)$.
	\end{enumerate}
	
	We present now the existence result based on asymptotic conditions in this second situation. 

\begin{theorem}
	\label{coro2.2}
Assume that conditions $(H_1)\text{--}(H_3)$ hold. Let $(k,l)\in\{(1,2),(2,1)\}$ and suppose that $\tilde{c}_k>0$, $\tilde{c}_l\geq0$, and either:
	\begin{enumerate}
		\item[\textup{$(a)$}] the $k$-th component satisfies $(\mathcal{A}_k)$ or $(\mathcal{B}_k)$, while the $l$-th component satisfies $(\widetilde{\mathcal{A}}_l)$;
		\item[\textup{$(b)$}] $M_k^l=M_l^k=0$, the $k$-th component satisfies $(\mathcal{A}_k)$ or $(\mathcal{B}_k)$, while the $l$-th component satisfies $(\widetilde{\mathcal{B}}_l)$.
	\end{enumerate}
Then, for every $\lambda_1,\lambda_2>0$, problem~\eqref{1.1} admits at least one positive solution $(u,v)$ such that $\|u\|_\infty>0$ and $\|v\|_\infty>0$.
\end{theorem}

\noindent
{\bf Proof.} Observe that the reasoning is similar to that in the proof of Theorem~\ref{coro2.1}, with changes only in the verification of condition \eqref{2.6}. Fix $\lambda_1,\lambda_2>0$ and let us prove $(a)$, assuming first that $\tilde c_1>0$, $\tilde c_2=0$ and that $(\mathcal{A}_1)$ and $(\widetilde{\mathcal{A}}_2)$ hold, with $m_j^i,M_j^i>0$ for $i,j\in\{1,2\}$. 

Since $(f_1+f_2)_0^1=+\infty$ uniformly with respect to $t\in I$ and $x_2\in\mathbb{R}^+$, there exists $\alpha_1>0$ for which condition~\eqref{2.5} is satisfied for the first variable, following the argument of the previous proof. Note that, since $\tilde c_1>0$ and $(f_1+f_2)_\infty^1=0$, uniformly with respect to $t\in I$ and $x_2\in\mathbb{R}^+$, we can choose a sufficiently large $\beta_1>\alpha_1$ such that
\[
\begin{aligned}
	(\lambda_1M_1^1f_1+\lambda_2M_2^1f_2)(t,x_1,x_2) &\leq \overline{\lambda}\,M^1(f_1+f_2)(t,x_1,x_2)\leq \frac{x_1}{B-A}, \\
	&\quad \text{for all } x_1\in[\tilde c_1 \beta_1,\infty), \, x_2\in\mathbb{R}_+ \text{ and } t\in I,
\end{aligned}
\]
which guarantees that \eqref{2.6} holds for $i=1$, and assumption $(\mathcal{C}_1)$ is not required. Moreover, following the reasoning in the proof of Theorem~\ref{coro2.1}, in order to verify the conditions for the second variable, it suffices to assume $(\widetilde{\mathcal{A}}_2)$ instead of $(\widetilde{\mathcal{A}}_2^0)$, since $\widetilde J_1=[\tilde c_1\alpha_1,\beta_1]$. This proves $(a)$ under $(\mathcal{A}_1)$ and $(\widetilde{\mathcal{A}}_2)$.

We next prove case~$(a)$ under $(\mathcal{B}_1)$ and $(\widetilde{\mathcal{A}}_2)$, with $\tilde{c}_1,\tilde{c}_2>0$. Since $(f_1+f_2)_0^1=0$ and $(f_1+f_2)_\infty^1=+\infty$ uniformly with respect to $t\in I$ and $x_2\in\mathbb{R}_+$, there exist constants $0<\alpha_1\equiv\alpha_1(\underline{\lambda})<\beta_1\equiv\beta_1(\overline{\lambda})$ such that
\[\overline{\lambda}\,M^1(f_1+f_2)(t,x_1,x_2)\leq \frac{x_1}{B-A}\leq \frac{\alpha_1}{B-A}, \, \text{ for all } x_1\in[0,\alpha_1], \, x_2\in\mathbb{R}_+ \text{ and } t\in I,
\]
and
\[\underline{\lambda}\,m^1(f_1+f_2)(t,x_1,x_2)\geq \frac{x_1}{c_1(b-a)}\geq \frac{\beta_1}{b-a}, \, \text{ for all } x_1\in[c_1\beta_1,\infty), \, x_2\in\mathbb{R}_+ \text{ and } t\in I.
\]
Consequently, \eqref{2.5} and \eqref{2.6} hold for the first variable.

Now, since $(f_1+f_2)_0^2=+\infty$ uniformly with respect to $t\in I$ and $x_1\in J_1:=[c_1\alpha_1,\beta_1]$, while $(f_1+f_2)_\infty^2=0$ uniformly with respect to $t\in I$ and $x_1\in \widetilde J_1:=[\tilde c_1 \alpha_1,\beta_1] $, there exist constants $0<\alpha_2\equiv\alpha_2(\underline{\lambda})<\beta_2\equiv\beta_2(\overline{\lambda})$ such that
\[\underline{\lambda}\,m^2(f_1+f_2)(t,x_1,x_2)\geq \frac{x_2}{c_2(b-a)}\geq \frac{\alpha_2}{b-a}, \, \text{ for all } x_1\in J_1, \, x_2\in(0,\alpha_2] \text{ and } t\in I,
\]
and
\[\overline{\lambda}\,M^2(f_1+f_2)(t,x_1,x_2)\leq \frac{x_2}{B-A}\leq \frac{\beta_2}{B-A}, \, \text{ for all } x_1\in\widetilde J_1, \, x_2\in[\tilde c_2\beta_2,\infty) \text{ and } t\in I.
\]
Therefore, conditions \eqref{2.5} and \eqref{2.6} hold for $i=2$. The remaining mixed configurations in the decoupled case, corresponding to $(b)$, follow by suitably combining the arguments above. \qed

\begin{remark}\label{rem:jde_compatibility}
	From a formal point of view, the compression--expansion conditions in Theorem~\ref{th1.1} could be considered for all four combinations of sublinear and superlinear behavior at zero and infinity that is, the combinations of conditions of type $\mathcal{A}$ and $\mathcal{B}$ for each variable. However, in Theorems~\ref{coro2.1} and \ref{coro2.2}, we restrict attention to the configurations that can actually be realized by continuous nonlinearities. The omitted coupled configurations would correspond to an empty class of admissible continuous nonnegative nonlinearities, since no such functions can satisfy the corresponding conditions simultaneously.
	
	To see this obstruction explicitly, consider $(\mathcal{B}_1)$ together with $(\widetilde{\mathcal{B}}_2)$ in the coupled case, that is, $M_2^1>0$ and hence $\delta_2^1=1$. Condition $(\mathcal{B}_1)$ requires that $(S_1^\delta)_0^1=0$ uniformly with respect to $(t,x_2)\in I\times \mathbb{R}^+$, namely,
	\[
	\lim_{x_1\to0^+}
	\sup_{(t,x_2)\in I\times\mathbb{R}^+}
	\frac{f_1(t,x_1,x_2)+f_2(t,x_1,x_2)}{x_1}
	=0.
	\]
	Taking $\varepsilon=1$, there exists $\delta_0>0$ such that $	(f_1+f_2)(t,x_1,x_2)<x_1$ for all $t\in I$, $x_1\in(0,\delta_0)$ and all $x_2\in\mathbb{R}^+$.
	
	Fix $x_1^*=\delta_0/2>0$. Then $J:=\{x_1^*\}$ is a compact subset of $\mathbb{R}_+$. For this fixed value of $x_1^*$, we have $(f_1+f_2)(t,x_1^*,x_2)\leq x_1^*$ for all $x_2\geq0$. Hence,
	\[
	\limsup_{x_2\to+\infty}
	\frac{S_2^\eta(t,x_1^*,x_2)}{x_2}
	\leq
	\limsup_{x_2\to+\infty}\frac{x_1^*}{x_2}
	=0.
	\]
	This contradicts condition $(\widetilde{\mathcal{B}}_2)$ (and likewise $(\widetilde{\mathcal{B}}_2^0)$), which requires $(S_2^\eta)_\infty^2=+\infty$ uniformly on every compact subset $J\subset\mathbb{R}_+$.

	By analogous arguments,  for $(k,l)\in \{(1,2), (2,1)\}$, we obtain the following:
	\begin{enumerate}[$(i)$]
		\item The configuration $(\mathcal{A}_k)$ together with $(\widetilde{\mathcal{B}}_l)$ (or $(\widetilde{\mathcal{B}}_l^0)$) is impossible when $M_l^k>0$. Indeed, the uniform limit at infinity in $(\mathcal{A}_k)$ bounds $S_k^\delta$ by $x_k$ for sufficiently large $x_k$, which precludes the superlinear growth of $S_l^\eta$ in the $x_l$-direction.
		
		\item The configuration $(\mathcal{B}_k)$ together with $(\widetilde{\mathcal{A}}_l^0)$ in Theorem~\ref{coro2.1} is impossible when $M_k^l>0$. The argument follows from the reasoning above.
	\end{enumerate}
	
	On the other hand, when the system is completely decoupled, that is, $M_l^k=M_k^l=0$, all off-diagonal weights vanish: $\eta_l^k=\delta_l^k=\eta_k^l=\delta_k^l=0$. Hence, the corresponding sums reduce to the individual nonlinearities. In this case, the growth conditions imposed on $f_k$ near zero and infinity place no restrictions on $f_l$. Consequently, all four pairings of conditions are compatible. In this sense, Theorem~\ref{coro2.1} extends results such as \cite[Corollary~4.1]{LFP_JRL}, which are recovered in the decoupled case, while allowing more general boundary conditions.
\end{remark}

Despite the mutual exclusivity of certain superlinear combinations in the coupled setting, Theorems~\ref{coro2.1} and \ref{coro2.2} extend and improve upon previous results \cite{Mou} by replacing global multivariable asymptotic assumptions, such as \cite[Hypothesis~(H8)]{Mou}, with the directional conditions $(S_k^\eta)_0^k$ and $(S_k^\delta)_\infty^k$. Moreover, it allows compressive and expansive behaviors to be combined along different directions. In addition, the solutions provided by Theorems~\ref{coro2.1} and \ref{coro2.2} are coexistence in the sense that both of their components are nontrivial, in contrast to those in, e.g., \cite{Mou,CIP,IMP,Y}, which may be semi-trivial.

On the other hand, while Theorem~\ref{coro2.1} requires directional asymptotic conditions that are uniform with respect to the cross-variable, coupled nonlinearities may instead be naturally described by their asymptotic behavior along the diagonal, that is, given a continuous function $g$, the following limits
\[(g)_0^{\mathrm{diag}}:=\lim_{r\to0^+}\frac{g(t,r,r)}{r},\qquad (g)_\infty^{\mathrm{diag}}:=\lim_{r\to+\infty}\frac{g(t,r,r)}{r}.\]	

When the nonlinearities are nondecreasing in each variable, these diagonal conditions yield the required estimates for existence, leading to the following result.

\begin{theorem}
	\label{th_diag}
	Assume that conditions $(H_1)$--$(H_3)$ hold. For each $i\in\{1,2\}$, suppose that $f_i(t,\cdot):\mathbb{R}^+\times \mathbb{R}^+\to\mathbb{R}^+$ is nondecreasing in each variable for every $t\in I$, and that
	\begin{equation}
		\label{3.28}
		(S_i^\eta)_0^{\mathrm{diag}}=+\infty
		\qquad\text{and}\qquad
		(S_i^\delta)_\infty^{\mathrm{diag}}=0,
		\quad\text{uniformly for }t\in I.
	\end{equation}
	Then, for every $\lambda_1,\lambda_2>0$, problem \eqref{1.1} admits at least one positive solution $(u,v)$ such that $	\|u\|_\infty>0$ and $\|v\|_\infty>0$.
\end{theorem}

\noindent
{\bf Proof.} Fix $\lambda_1,\lambda_2>0$. We restrict the proof to the case
$m_j^i>0$ for $i,j\in\{1,2\}$, since the remaining cases are
analogous. For each $i\in\{1,2\}$, the condition $(S_i^\eta)_0^{\mathrm{diag}}=+\infty$, uniformly for $t\in I$, yields the existence of $\alpha\equiv\alpha(\lambda_1,\lambda_2)>0$ sufficiently small such that
\[(f_1+f_2)(t,r,r)>\frac{r}{\underline{\lambda}\,m^ic(b-a)},\qquad \text{ for all } r\in(0,\alpha] \text{ and } t\in I,\quad i\in\{1,2\},\]
where $c=\min\{c_1,c_2\}$.
The monotonicity of $f_1$ and $f_2$ with respect to both variables then gives
\[\begin{aligned}
(\lambda_1m_1^i f_1+\lambda_2m_2^i f_2)(t,x_1,x_2) &\geq \underline{\lambda}\,m^i(f_1+f_2)(t,x_1,x_2)\\
&\geq \underline{\lambda}\,m^i(f_1+f_2)(t,c_1 \alpha,c_2 \alpha)\\
&\geq\underline{\lambda}\,m^i(f_1+f_2)(t,c\alpha,c\alpha)\\
&>\frac{\alpha}{b-a}, \quad \text{ for all } (x_1,x_2)\in[c_1\alpha,\infty)\times [c_2\alpha,\infty) \text{ and } t\in I, \quad i\in\{1,2\}.
\end{aligned}\]

On the other hand, the condition
$(S_i^\delta)_\infty^{\mathrm{diag}}=0$, uniformly for $t\in I$, for
each $i\in\{1,2\}$, allows us to choose
$\beta\equiv\beta(\lambda_1,\lambda_2)>0$ sufficiently large ($\beta>\alpha$) so that
\[(f_1+f_2)(t,r,r)<\frac{r}{\overline{\lambda}\,M^i(B-A)},\qquad \text{ for all }r\in[\beta,\infty) \text{ and } t\in I,\quad i\in\{1,2\}.
\]
Thus, the monotonicity of the nonlinearities implies
\[\begin{aligned}
(\lambda_1M_1^i f_1+\lambda_2M_2^i f_2)(t,x_1,x_2)&\leq\overline{\lambda}\,M^i (f_1+f_2)(t,x_1,x_2)\\
&\leq\overline{\lambda}\,M^i(f_1+f_2)(t,\beta,\beta)\\
&<\frac{\beta}{B-A}, \quad \text{ for all } (x_1,x_2)\in[0,\beta]^2 \text{ and } t\in I, \quad i\in\{1,2\}. 
\end{aligned}\]

Hence, conditions \eqref{2.5} and \eqref{2.6} of Theorem~\ref{th2.1} are fulfilled by taking $r_i=\alpha$ and $R_i=\beta$ for $i\in\{1,2\}$, and the desired conclusion follows. \qed

The following example illustrates the different configurations covered by Theorems~\ref{coro2.1} and \ref{coro2.2}, as well as a setting in which Theorem~\ref{th_diag} applies while Theorem~\ref{coro2.1} and \ref{coro2.2} does not.
	\begin{example}
		\label{ex4.1}
	Let $0<A<B<\infty$ be two fixed radios, and set $L:=B-A$. In this example, we provide several concrete choices of nonlinearities $f_1,f_2\colon\mathbb{R}^+\times\mathbb{R}^+\to\mathbb{R}^+$ for which the elliptic system posed on the annular domain $\Omega=\{x\in\mathbb{R}^3:A<|x|<B\}$ with unit width $L=1$ given by
		\begin{equation}\label{eq:pde_system2}
			\left\{
			\begin{aligned}
				-\Delta U &= \lambda_1  f_1\big(U, V\big), & x\in \Omega, \\[0.6em]
				-\Delta V &= \lambda_2  f_2\big(U, V\big), & x\in \Omega,
			\end{aligned}
			\right.
		\end{equation}
		 and subject to the coupled Robin boundary conditions
		\begin{equation}\label{eq:pde_bc12}
			\begin{aligned}
				U(x) - \frac{A}{B} \frac{\partial U}{\partial \nu}(x) + \alpha V(x) &= 0, \quad |x| = A; \qquad & U(x) + \frac{B}{A} \frac{\partial U}{\partial \nu}(x) &= 0, \quad |x| = B, \\[0.8em]
				V(x) - \beta \frac{A}{B} \frac{\partial V}{\partial \nu}(x) + \alpha U(x) &= 0, \quad |x| = A; \qquad & V(x) + \beta \frac{B}{A} \frac{\partial V}{\partial \nu}(x) &= 0, \quad |x| =B,
			\end{aligned}
		\end{equation}
	where $\alpha < 0$, $\beta \ge 0$, and $D(\alpha,\beta) := 2\alpha^2(\beta+1) - 1 > 0$, admits at least one radially symmetric solution $(U, V)$---that is, $U(x) = U(|x|)$ and $V(x) = V(|x|)$ for all $x \in \overline{\Omega}$---having both components nontrivial.
	
		
The radially symmetric solutions of this problem, with $r=|x|$ and after the change of variables $r:I\to I$, $I:=[A,B]$, given by
\[r(t)=\dfrac{AB}{A+B-t}, \]
correspond to solutions of the following system of ordinary differential equations
\begin{equation}\label{eq:ode_system2}
				\left\{
	\begin{aligned}
		-u''(t) &= \lambda_1 h(t) f_1\big(u(t), v(t)\big), \quad t \in I, \\[0.6em]
		-v''(t) &= \lambda_2 h(t) f_2\big(u(t), v(t)\big), \quad t \in I, \\[0.6em]
		u(A) + u'(A) + &\alpha v(A) = 0, \qquad u(B) + u'(B) = 0, \\[0.6em]
		v(A) + \beta v'(A) + &\alpha u(A) = 0, \qquad v(B) + \beta v'(B) = 0,
	\end{aligned}
				\right.
\end{equation}
where $h:[A,B]\to\mathbb{R}$ is the weight function defined as
\[h(t):=\dfrac{A^2 B^2}{(A+B-t)^4}. \]
Observe that $0<A^2/B^2\leq h(t)\leq B^2/A^2$ for all $t\in I$, so $h$ plays no a relevant role in the verification of the asymptotic conditions required by Theorems~\ref{coro2.1}, \ref{coro2.2} or \ref{th_diag}.

In this case, the four Green's functions are given by
\begin{equation*}
	\begin{aligned}
		G_1^1(t, s) &= \frac{\big(B + 1 - \max\{r, s\}\big)\big[1 - C_1(\alpha, \beta)(\min\{r, s\} - A)\big]}{D(\alpha, \beta)}, && (t, s) \in I \times I, \\[0.6em]
		G^1_2(t, s) &= \frac{-\alpha \beta (B + 1 - r)(B + \beta - s)}{D(\alpha, \beta)}, && (t, s) \in I \times I, \\[0.6em]
		G^2_1(t, s) &= \frac{-\alpha (B + \beta - r)(B + 1 - s)}{D(\alpha, \beta)}, && (t, s) \in I \times I, \\[0.6em]
		G^2_2(t, s) &= \frac{\big(B + \beta - \max\{r, s\}\big)\big[\beta - C_2(\alpha, \beta)(\min\{r, s\} - A)\big]}{D(\alpha, \beta)}, && (t, s) \in I \times I,
	\end{aligned}
\end{equation*}
where $C_1(\alpha,\beta):=1-(1+\beta)\alpha^2$ and $C_2(\alpha,\beta):= 1-2\alpha^2$.

Under the prescribed assumptions on the parameters $\alpha$ and $\beta$, the four Green's functions satisfy hypotheses $(H_1)$ and $(H_2)$; see Appendices~\ref{app:A}, \ref{app:B}, and~\ref{app:C} for details. Moreover, if $\beta>0$, we can take $[a,b]=[A,B]$, yielding
$c_1=\tilde{c}_1>0$ and $c_2=\tilde{c}_2>0$. If $\beta=0$, however, the presence of a Dirichlet boundary condition requires us to take a proper subinterval $[a,b]\subsetneq[A,B]$, in which case $c_1,\tilde{c}_1,c_2>0$ and $\tilde{c}_2=0$.

\begin{enumerate}[$(i)$]
	\item Under the assumptions $\alpha < 0$, $\beta \ge 0$, and $D(\alpha,\beta) > 0$, problem~\eqref{eq:ode_system2} admits at least one solution with both components nontrivial for any $\lambda_1, \lambda_2 > 0$, whenever the nonlinearities are given by
	\[
	f_1(x_1, x_2) = \frac{x_1^{p_1}}{1 + x_1^{q_1}} \left( 1 + \frac{x_2^{\alpha_1}}{1 + x_2^{\alpha_1}} \right), 
	\qquad
	f_2(x_1, x_2) = \frac{x_2^{p_2}}{1 + x_2^{q_2}} \left( 1 + \frac{x_1^{\alpha_2}}{1 + x_1^{\alpha_2}} \right),
	\]
	with $0 < p_1, p_2 < 1$, $q_i \ge p_i$, and $\alpha_i > 0$ for $i=1,2$.
	
First, note that $0 \le f_1(x_1, x_2) + f_2(x_1, x_2) \le 4$ for all $(x_1, x_2) \in \mathbb{R}^+\times \mathbb{R}^+$, so condition~$(\mathcal{C}_1)$ is clearly fulfilled.
Next, since $1 - p_1 > 0$, we have
\[
\inf_{x_2 \ge 0} \frac{(f_1 + f_2)(x_1, x_2)}{x_1} 
\ge \frac{1}{x_1^{1 - p_1}(1 + x_1^{q_1})} \xrightarrow[x_1 \to 0^+]{} +\infty,
\]
which implies that $(S_1^\eta)_0^1 = +\infty$ uniformly on $\mathbb{R}^+$. On the other hand,
\[
\sup_{x_2 \ge 0} \frac{(f_1 + f_2)(x_1, x_2)}{x_1} \le \frac{4}{x_1} \xrightarrow[x_1 \to +\infty]{} 0,
\]
yielding $(S_1^\delta)_\infty^1 = 0$ uniformly on $\mathbb{R}^+$. Thus, condition~$(\mathcal{A}_1)$ holds.
Similarly, let $J \subset \mathbb{R}^+$ be an arbitrary compact set. Since $1 - p_2 > 0$, we obtain
\[
\inf_{x_1 \in J} \frac{(f_1 + f_2)(x_1, x_2)}{x_2} \ge \frac{1}{x_2^{1 - p_2}(1 + x_2^{q_2})} \xrightarrow[x_2 \to 0^+]{} +\infty,
\, \text{ and } \,
\sup_{x_1 \in J} \frac{(f_1 + f_2)(x_1, x_2)}{x_2} \le \frac{4}{x_2} \xrightarrow[x_2 \to +\infty]{} 0,
\]
from which $(\widetilde{\mathcal{A}}_2^0)$ follows, and consequently $(\widetilde{\mathcal{A}}_2)$ is satisfied. Therefore, Theorems~\ref{coro2.1} or \ref{coro2.2} can be applied.

\item Under the assumptions $\alpha < 0$, $\beta > 0$, and $D(\alpha,\beta) > 0$, problem~\eqref{eq:ode_system2} admits at least one solution with both components nontrivial for any $\lambda_1, \lambda_2 > 0$, whenever the nonlinearities are given by
\[
f_1(x_1, x_2) = f_2(x_1, x_2) = x_1^p \bigl( 1 + \sin^2(x_2) \bigr),
\]
with $p > 1$.

First, since $0 \le \sin^2(x_2) \le 1$ for all $x_2 \in \mathbb{R}^+$, we have $2 x_1^p \le (f_1 + f_2)(x_1, x_2) = 2 x_1^p \bigl( 1 + \sin^2(x_2) \bigr) \le 4 x_1^p$ for all $(x_1, x_2) \in \mathbb{R}^+\times \mathbb{R}^+$. Moreover, since $p - 1 > 0$, we have
\[
\sup_{x_2 \ge 0} \frac{(f_1 + f_2)(x_1, x_2)}{x_1} \le 4 x_1^{p - 1} \xrightarrow[x_1 \to 0^+]{} 0,
\]
which implies that $(S_1^\delta)_0^1 = 0$ uniformly on $\mathbb{R}^+$. On the other hand,
\[
\inf_{x_2 \ge 0} \frac{(f_1 + f_2)(x_1, x_2)}{x_1} \ge 2 x_1^{p - 1} \xrightarrow[x_1 \to +\infty]{} +\infty,
\]
yielding $(S_1^\eta)_\infty^1 = +\infty$ uniformly on $\mathbb{R}^+$. Thus, hypothesis~$(\mathcal{B}_1)$ holds.

Similarly, let $J \subset \mathbb{R}_+$ be an arbitrary compact set, and choose constants $0 < a_J \le b_J < \infty$ such that $a_J \le x_1 \le b_J$ for all $x_1 \in J$. Since $a_J > 0$ and $1 + \sin^2(x_2) \ge 1$,
\[
\inf_{x_1 \in J} \frac{(f_1 + f_2)(x_1, x_2)}{x_2} \ge \frac{2 a_J^p}{x_2} \xrightarrow[x_2 \to 0^+]{} +\infty,
\qquad
\sup_{x_1 \in J} \frac{(f_1 + f_2)(x_1, x_2)}{x_2} \le \frac{4 b_J^p}{x_2} \xrightarrow[x_2 \to +\infty]{} 0,
\]
from which $(\widetilde{\mathcal{A}}_2)$ holds. Then Theorem~\ref{coro2.2} gives the conclusion again.

\item  Under the assumptions $\alpha < 0$, $\beta \ge 0$, and $D(\alpha,\beta) > 0$, problem~\eqref{eq:ode_system2} admits at least one solution with both components nontrivial for any $\lambda_1, \lambda_2 > 0$, whenever the nonlinearities are given by
\[f_1(x_1, x_2) =x_1^{p_1} x_2^{q_1}, \qquad  f_2( x_1, x_2) = x_1^{p_2} x_2^{q_2},\]
with $p_i, q_i>0$ and $p_i+q_i<1$, $i=1,2$.

Clearly, $f_1$ and $f_2$ are nondecreasing in both variables, as $p_i, q_i>0$, $i=1,2$. Moreover, for every $r>0$,
\[\frac{(f_1+f_2)(r,r)}{r}=r^{p_1+q_1-1}+r^{p_2+q_2-1}.\]
The assumption $p_i+q_i<1$, $i=1,2$, then yields the limits in \eqref{3.28} as $r\to0$ and $r\to\infty$. Therefore, Theorem~\ref{th_diag} gives the conclusion.
\end{enumerate}
Observe that, in this example, unlike in the application of Theorem~\ref{th2.1} in Example~\ref{ej:2.1}, Theorems~\ref{coro2.1}, \ref{coro2.2} and \ref{th_diag} allow us to establish existence without fixing $\alpha$ and $\beta$ to obtain bounds from the extrema of the Green's functions. It suffices to determine whether $\tilde{c}_1$ and $\tilde{c}_2$ are positive or vanish. Note that, in case~(ii), we therefore require $\beta>0$ in order to apply Theorem~\ref{coro2.2}.
\end{example}

\section{A multiplicity result}
\label{sec:multiplicity}

The existence results established in the preceding sections guarantee the existence of at least one positive solution. In Section~4, this is achieved either through directional asymptotic conditions (Theorems~\ref{coro2.1} and~\ref{coro2.2}) or through diagonal growth conditions combined with monotonicity (Theorem~\ref{th_diag}). A natural question is whether these two types of behavior can be combined to establish the existence of multiple positive solutions. In this section, we answer this question affirmatively by establishing a multiplicity result that guarantees the existence of at least two positive solutions, both of whose components are nontrivial. 

\begin{theorem}
	\label{th:multiplicity_two_solutions}
	Assume that conditions $(H_1)$--$(H_3)$ hold. Suppose, moreover, that for each $k\in\{1,2\}$,
	\begin{equation}
		\label{eq:mult_diag_zero}
		(S_k^\eta)_0^{\mathrm{diag}} = +\infty,
		\quad \text{uniformly with respect to } t\in I,
	\end{equation}
	and that the following conditions are satisfied:
	\begin{enumerate}
		\item[(i)] For every $t\in I$, the function
		$f_k(t,\cdot,\cdot):[0,\infty)^2\to[0,\infty)$
		is nondecreasing in each variable.
		
		\item[(ii)] There exists $r^*>0$ such that
		\[
		(S_k^\eta)_\infty^k=+\infty
		\quad \text{for all } x_l\in(0,r^*],\ t\in I,\ k\neq l.
		\]
	\end{enumerate}
	Then, for all sufficiently small $\lambda_1,\lambda_2>0$, problem~\eqref{1.1} admits at least two positive solutions, each having both components nontrivial.
\end{theorem}

\noindent
\textbf{Proof.}  For simplicity, we present the proof under the assumption that $m_j^i>0$ for all $i,j\in\{1,2\}$, since the other cases can be treated analogously. We begin by proving that there exist $r, R \in\mathbb{R} $ with $0 < r < r^* <R$ such that $\|\mathcal{T}_i x\|_\infty > \|x_i\|_\infty$ whenever $\|x_i\|_\infty \in \{r, R\}$ for $i \in \{1, 2\}$. Indeed, by \eqref{eq:mult_diag_zero} and the monotonicity of $f_1$ and $f_2$ with respect to $x_1$ and $x_2$, for each pair $\lambda_1, \lambda_2 > 0$ we can choose $r \equiv r(\lambda_1, \lambda_2) > 0$ sufficiently small such that $r < r^*$ and
\begin{equation}\label{eq:mult_exp_r}
(\lambda_1 m_1^i f_1 + \lambda_2 m_2^i f_2)(t, x_1, x_2) \ge \underline{\lambda} m^i (f_1 + f_2)(t, cr, cr) > \frac{r}{b-a}, \text{ for all } (x_1, x_2) \in [c r, \infty)^2 \text{ and } t \in I, \, i \in \{1, 2\},
\end{equation}
where $c := \min\{c_1, c_2\}$.

On the other hand, the monotonicity of the nonlinearities combined with condition~$(ii)$ implies the existence of a sufficiently large radius $R> 0$ such that
\begin{equation}\label{eq:mult_exp_R}
	\begin{aligned}
	(\lambda_1 m_1^i f_1 + \lambda_2 m_2^i f_2)(t, x_1, x_2) &\ge \underline{\lambda} m^i (f_1 + f_2)(t, x_1, x_2) \\
	&> \frac{R}{b-a}, \text{ for all } x_i \in [c_i R, \infty), \, x_j \in [c_j r, \infty), \, j\ne i \text{ and } i \in \{1, 2\}.
	\end{aligned}
\end{equation} 

Therefore, in view of the estimates in Theorem~\ref{th2.1}, for each $(x_1, x_2) \in \overline{K}_{(r,r), (R,R)}$ we have that
\[
\|\mathcal{T}_i x\|_\infty > \|x_i\|_\infty, \quad \text{if } \|x_i\|_\infty \in \{r, R\}, \quad i \in \{1, 2\},
\]
which implies $i_{K}(\mathcal{T}, K_{(r,r), (R,R)}) = 0$. Indeed, let $\mathcal{N} : \overline{K}_{(R,R)} \to K$ be an extension of $\mathcal{T}$ of the form $\mathcal{N} = \mathcal{T} \circ \rho$, where $\rho = (\rho_1, \rho_2)$ and each $\rho_i : (\overline{K}_i)_{R} \to (\overline{K}_i)_{r, R}$ is a retraction (see \cite{LFP_JRL}). As thus defined, for each $(x_1, x_2) \in \overline{K}_{(R,R)}$ it holds that
\[
\|\mathcal{N}_i x\|_\infty > \|x_i\|_\infty, \quad \text{if } \|x_i\|_\infty \in \{r, R\}, \quad i \in \{1, 2\}.
\]
In this way, by Proposition \ref{prop2.2} we have that
\[
i_K(\mathcal{N}, K_{(r,r)}) = i_K(\mathcal{N}, K_{(R,R)}) = i_K\bigl(\mathcal{N}, (K_1)_{r} \times (K_2)_{R}\bigr) = i_K\bigl(\mathcal{N}, (K_1)_{R} \times (K_2)_{r}\bigr) = 0,
\]
and thus, using the additivity property of the fixed point index together with the fact that $\mathcal{N} = \mathcal{T}$ on $K_{(r,r), (R,R)}$, we obtain
\[
\begin{aligned}
i_K(\mathcal{T}, K_{(r,r), (R,R)}) = i_K(\mathcal{N}, K_{(r,r), (R,R)}) =& i_K(\mathcal{N}, K_{(R,R)}) - i_K\bigl(\mathcal{N}, (K_1)_{r} \times (K_2)_{R}\bigr) \\
&- \Bigl[i_K\bigl(\mathcal{N}, (K_1)_{R} \times (K_2)_{r}\bigr) - i_K(\mathcal{N}, K_{(r,r)})\Bigr] = 0.
\end{aligned}
\]

Next, we choose $\lambda_1, \lambda_2 > 0$ sufficiently small such that
\[
\lambda_1 M_1^i f_1(t, r^*, r^*) + \lambda_2 M_2^i f_2(t, r^*, r^*) < \frac{r^*}{B-A}, \quad \text{for all } t \in I, \quad  i \in \{1, 2\}.
\]
By the monotonicity of $f_1$ and $f_2$ with respect to $x_1$ and $x_2$, this choice implies that
\[
\begin{aligned}
(\lambda_1 M_1^i f_1 + \lambda_2 M_2^i f_2)(t, x_1, x_2) &\le (\lambda_1 M_1^i f_1 + \lambda_2 M_2^i f_2)(t, r^*, r^*) \\
&< \frac{r^*}{B-A}, \text{ for all } (x_1, x_2) \in [0, r^*]^2 \text{ and } t \in I, \quad i \in \{1, 2\}.
\end{aligned}
\]
Consequently, for each $(x_1, x_2) \in \overline{K}_{(r,r), (r^*,r^*)}$ we have
\[
\|\mathcal{T}_i x\|_\infty > \|x_i\|_\infty \quad \text{if } \|x_i\|_\infty = r, \qquad \text{and} \qquad \|\mathcal{T}_i x\|_\infty < \|x_i\|_\infty \quad \text{if } \|x_i\|_\infty = r^*,
\]
from which it follows directly that $i_K(\mathcal{T}, K_{(r,r), (r^*,r^*)}) = 1$ (see Remark \ref{rem2.1}).

\begin{figure}[h]
	\centering
	\includegraphics[scale=0.9]{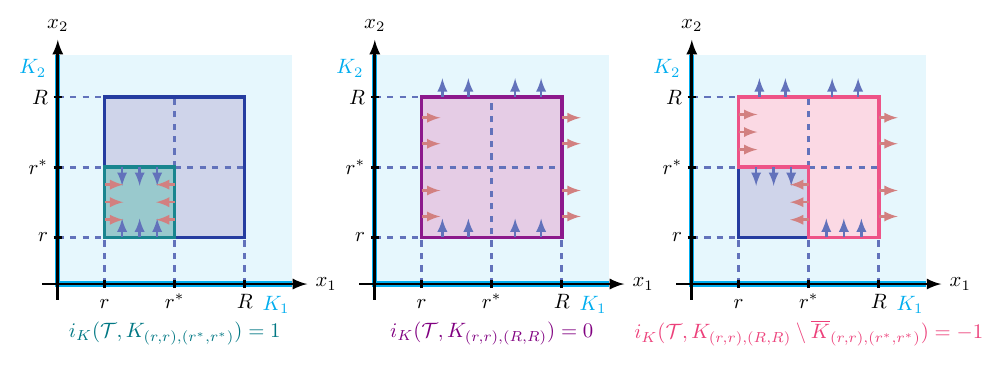}
	\vspace{-0.2cm}
	\caption{Geometric visualization of the fixed point index computations in Theorem \ref{th:multiplicity_two_solutions}. }
	\label{fig:index_computation}
\end{figure}

Finally, applying once again the additivity property of the fixed point index, we find that
\[
i_K(\mathcal{T}, K_{(r,r), (R,R)} \setminus \overline{K}_{(r,r), (r^*,r^*)}) = i_K(\mathcal{T}, K_{(r,r), (R,R)}) - i_K(\mathcal{T}, K_{(r,r), (r^*,r^*)}) = 0 - 1 = -1.
\]
See Figure~\ref{fig:index_computation} for a geometric illustration of the fixed point index computations on the corresponding sets.

By the existence property of the fixed point index, we conclude that problem~\eqref{1.1} admits at least two positive solutions with both components nontrivial, one located in $K_{(r,r), (r^*,r^*)}$ and the other in $K_{(r,r), (R,R)} \setminus \overline{K}_{(r,r), (r^*,r^*)}$.
\qed

We conclude this section by illustrating the applicability of Theorem~\ref{th:multiplicity_two_solutions} to a coupled system of ordinary differential equations with two parameters and power type nonlinearities.

\begin{example}\label{ej:5.1}
The system of ordinary differential equations
\begin{equation}\label{eq:ode_mult_example}
\left\{
\begin{aligned}
	-u''(t) &= \lambda_1 \left(u^{p_1}(t) + v^{q_1}(t)\right), \quad t \in I := [1, 2], \\[0.1em]
	-v''(t) &= \lambda_2 \left(u^{q_2}(t) + v^{p_2}(t)\right), \quad t \in I, \\[0.1em]
	u(1) + u'(1) - v(1) &= 0, \qquad u(2) + u'(2) = 0, \\[0.1em]
	v(1) + v'(1) - u(1) &= 0, \qquad v(2) + v'(2) = 0,
\end{aligned}
\right.
\end{equation}
with $0<p_i<1$ and $q_i>1$, $i\in\{1,2\}$, admits at least two distinct positive solutions with both components nontrivial whenever
\begin{equation}\label{eq:region_M}
(\lambda_1, \lambda_2) \in \mathcal{M} := \left\{ (x_1, x_2) \in \mathbb{R}^2 : x_1 > 0, \ x_2 > 0, \ x_2 < \frac{3}{8} - \frac{9}{16} x_1, \text{ and } x_2 < \frac{2}{3} - \frac{16}{9} x_1 \right\},
\end{equation}
as a consequence of Theorem~\ref{th:multiplicity_two_solutions}.

Note that system~\eqref{eq:ode_mult_example} is of the form of system~\eqref{eq:ode_system} considered in Example~\ref{ej:2.1}, but with two independent parameters and with nonlinearities
\[
f_1(x_1, x_2) := x_1^{p_1} + x_2^{q_1} \qquad \text{and} \qquad f_2(x_1, x_2) := x_1^{q_2} + x_2^{p_2}.
\]
Indeed, taking $[a, b] = I$, the corresponding Green's functions $G_j^i$ are continuous and strictly positive on $I^2$ and satisfy conditions $(H_1)$ and $(H_2)$ with $c_j = \tilde{c}_j = 1/2$, $m_j^i = 1/3$, $M_2^1 = M_1^2 = 4/3$, and $M_1^1 = M_2^2 = 3/4$ for $i, j \in \{1, 2\}$.

The nonlinearities $f_1$ and $f_2$ are clearly nondecreasing in both variables. Moreover,
\[
(f_1 + f_2)(x_1, x_2) = x_1^{p_1} + x_2^{p_2} + x_1^{q_2} + x_2^{q_1}, \quad (x_1, x_2) \in \mathbb{R}^+ \times \mathbb{R}^+,
\]
and, for every $r > 0$,
\[
\frac{(f_1 + f_2)(r, r)}{r} = r^{p_1-1} + r^{p_2-1} + r^{q_2-1} + r^{q_1-1}.
\]
Consequently, since $p_i-1<0$, $\lim_{r \to 0^+} \frac{(f_1 + f_2)(r, r)}{r} = +\infty$, and \eqref{eq:mult_diag_zero} follows.

Fixing now $r^* = 1 > 0$, observe that
\[
(f_1 + f_2)_\infty^1 = \lim_{x_1 \to \infty} \left( \frac{1}{x_1^{1-p_1}} + \frac{x_2^{p_2}}{x_1} + x_1^{q_2-1} + \frac{x_2^{q_1}}{x_1} \right) \ge \lim_{x_1 \to \infty}  x_1^{q_2-1} = +\infty,
\]
whence $(f_1 + f_2)_\infty^1 = +\infty$ for $x_2 \in [0, 1]$. By a symmetric reasoning, it follows that $(f_1 + f_2)_\infty^2 = +\infty$ for $x_1 \in [0, 1]$. 

Following the proof of Theorem~\ref{th:multiplicity_two_solutions}, and given the choice $r^* = 1$, the upper compressive condition at $r^*$ requires
\[
\lambda_1 M_1^i f_1(1, 1) + \lambda_2 M_2^i f_2(1, 1) < \frac{r^*}{B - A} = 1, \quad i \in \{1, 2\}.
\]
Since $f_1(1, 1) = f_2(1, 1) = 2$, we conclude the existence of the two desired positive solutions whenever $(\lambda_1, \lambda_2) \in \mathcal{M}$.
\end{example}

\section*{Acknowledgements}

L. M. Fern\'andez-Pardo acknowledges the financial support of the Spanish Ministry of Science, Innovation and Universities (Grant reference FPU24/01545). Both authors have been partially supported by Grant PID2025-169849NB-I00 funded by MICIU/AEI/10.13039/501100011033, and by Project ED431C 2023/12 funded by Xunta de Galicia, Spain.

\appendix

\section{Explicit derivation of the Green's matrix and kernels in Example \ref{ex4.1}}
\label{app:A}
Consider the coupled second-order boundary value problem \eqref{eq:ode_system2} with parameters $\alpha, \beta \in \mathbb{R}$, constants $\lambda_1, \lambda_2 > 0$, and unit interval length $B- A = 1$.

Introducing the state vector $z(t):=(u(t), u'(t), v(t), v'(t))^T \in \mathbb{R}^4$, system  \eqref{eq:ode_system2} transforms into:
\begin{equation}\label{eq:first_order}
	z'(t) = \mathcal{A} z(t) + \mathcal{F}(t), \quad t \in I, \qquad E z(A) + F z(B) = 0,
\end{equation}
where the system matrix $\mathcal{A} \in \mathcal{M}_{4\times 4}(\mathbb{R})$ and source vector $\mathcal{F}(t) \in \mathbb{R}^4$ are:
\begin{equation}\label{eq:A_f}
	\mathcal{A} = \begin{pmatrix}
		0 &\quad 1 & \quad 0 & \quad 0 \\
		0 &\quad  0 &\quad  0 &\quad  0 \\
		0 &\quad  0 &\quad  0 &\quad  1 \\
		0 &\quad  0 &\quad  0 &\quad  0
	\end{pmatrix}, \qquad
	\mathcal{F}(r) = \begin{pmatrix}
		0 \\
		-\lambda_1 h(t) f_1(u(t), v(t)) \\
		0 \\
		-\lambda_2 h(t) f_2(u(t), v(t))
	\end{pmatrix}.
\end{equation}
According to \cite{Mou}, the boundary conditions in \eqref{eq:first_order} are encoded by the matrices $E, F \in \mathcal{M}_{4\times 4}(\mathbb{R})$ given by
\begin{equation}\label{eq:E_F}
	E = \begin{pmatrix}
		1 &\quad  1 &\quad  \alpha &\quad  0 \\
		0 &\quad  0 &\quad  0 &\quad  0 \\
		\alpha &\quad  0 &\quad  1 &\quad  \beta \\
		0 &\quad  0 &\quad  0 &\quad  0
	\end{pmatrix}, \qquad
	F = \begin{pmatrix}
		0 &\quad  0 &\quad  0 &\quad  0 \\
		1 &\quad  1 &\quad  0 &\quad  0 \\
		0 &\quad  0 &\quad  0 &\quad  0 \\
		0 &\quad  0 &\quad  1 &\quad  \beta
	\end{pmatrix}.
\end{equation}
Since $\mathcal{A}^2 = 0$, the fundamental matrix $\Phi(t) = e^{\mathcal{A}t} = I_4 + \mathcal{A}t$ and its inverse $\Phi^{-1}(s)$ are
\begin{equation}\label{eq:Phi}
	\Phi(t) = \begin{pmatrix}
		1 &\quad  t &\quad  0 &\quad  0 \\
		0 &\quad  1 &\quad  0 &\quad  0 \\
		0 &\quad  0 &\quad  1 &\quad  t \\
		0 &\quad  0 &\quad  0 &\quad  1
	\end{pmatrix}, \qquad
	\Phi^{-1}(s) = \begin{pmatrix}
		1 &\quad  -s &\quad  0 &\quad  0 \\
		0 &\quad  1 &\quad  0 &\quad  0 \\
		0 &\quad  0 &\quad  1 &\quad  -s \\
		0 &\quad  0 &\quad  0 & \quad 1
	\end{pmatrix}.
\end{equation}

\subsection{Characteristic matrix $M_\Phi$, determinant, and inverse $(M_\Phi)^{-1}$}

Following \cite{Mou}, the boundary characteristic matrix $M_\Phi:= E \Phi(A) + F \Phi(B)$ is
\begin{equation}\label{eq:MPhi}
	M_\Phi = \begin{pmatrix}
		1 &\quad  A + 1 &\quad  \alpha & \quad \alpha A \\
		1 &\quad  B + 1 &\quad  0 &\quad  0 \\
		\alpha &\quad  \alpha A &\quad  1 &\quad  A + \beta \\
		0 &\quad  0 &\quad  1 &\quad  B + \beta
	\end{pmatrix}.
\end{equation}
Under the condition $B - A = 1$, the determinant simplifies to
\begin{equation}\label{eq:det}
	\det(M_\Phi) = -(2\alpha^2\beta + 2\alpha^2 - 1) = -D(\alpha, \beta), \qquad \text{with} \quad D\equiv D(\alpha, \beta) := 2\alpha^2(1 + \beta) - 1.
\end{equation}
The homogeneous problem has only the trivial solution if and only if $D(\alpha, \beta) \neq 0$.

Introducing the auxiliary parameters:
\begin{equation}\label{eq:C1_C2}
	C_1\equiv C_1(\alpha, \beta) :=1 - \alpha^2(1 + \beta) = \frac{1 + D(\alpha, \beta)}{2}, \qquad C_2\equiv C_2(\alpha, \beta) := 1 - 2\alpha^2,
\end{equation}
the inverse $(M_\Phi)^{-1} = -\frac{1}{D(\alpha, \beta)} \operatorname{adj}(M_\Phi)$ is given explicitly by
\begin{equation}\label{eq:MPhi_inv}
	(M_\Phi)^{-1} = \frac{1}{D(\alpha, \beta)} \begin{pmatrix}
		-(A + 2) &\quad  1 + A C_1 & \alpha(A + 2)(\beta + 1) &\quad  -\alpha\beta(A + 2) \\[0.25em]
		1 &\quad  -C_1 & -\alpha(\beta + 1) &\quad  \alpha\beta \\[0.25em]
		2\alpha(A + \beta + 1) &\quad  -\alpha(A + \beta + 1) & -(A + \beta + 1) & \quad  \beta + A C_2 \\[0.25em]
		-2\alpha &\quad  \alpha & 1 & \quad  -C_2
	\end{pmatrix}.
\end{equation}

The boundary projection matrices $M_E := (M_\Phi)^{-1} E \Phi(A)$ and $M_F :=(M_\Phi)^{-1} F \Phi(B)$ satisfy
\begin{equation}\label{eq:partition_I}
	M_E + M_F = (M_\Phi)^{-1} [E \Phi(A) + F \Phi(B)] = (M_\Phi)^{-1} M_\Phi = I_4.
\end{equation}
Carrying out the products and simplifying with $B = A + 1$, we have
\begin{equation}\label{eq:ME}
	M_E = \frac{1}{D(\alpha, \beta)} \begin{pmatrix}
		-(B + 1) C_1 &\quad  -(B + 1)(1 + A C_1) &\quad  \alpha\beta(B + 1) &\quad  \alpha\beta(B + 1)(B + \beta) \\[0.25em]
		C_1 &\quad  1 + A C_1 &\quad  -\alpha\beta &\quad  -\alpha\beta(B + \beta) \\[0.25em]
		\alpha(B + \beta) &\quad  \alpha(B + 1)(B + \beta) &\quad  -(B + \beta)C_2 &\quad  -(B + \beta)(\beta + A C_2) \\[0.25em]
		-\alpha &\quad  -\alpha(B + 1) &\quad  C_2 &\quad  \beta + A C_2
	\end{pmatrix}.
\end{equation}

Complementarily, via $M_F = I_4 - M_E$, we obtain
\begin{equation}\label{eq:MF}
	M_F = \frac{1}{D(\alpha, \beta)} \begin{pmatrix}
		1 + A C_1 &\quad  (B + 1)(1 + A C_1) &\quad  -\alpha\beta(B + 1) &\quad  -\alpha\beta(B + 1)(B + \beta) \\[0.25em]
		-C_1 &\quad  -(B + 1) C_1 &\quad  \alpha\beta &\quad  \alpha\beta(B + \beta) \\[0.25em]
		-\alpha(B + \beta) &\quad  -\alpha(B + 1)(B + \beta) &\quad  \beta + A C_2 & (B + \beta)(\beta + A C_2) \\[0.25em]
		\alpha &\quad  \alpha(B + 1) &\quad  -C_2 &\quad  -(B + \beta) C_2
	\end{pmatrix}.
\end{equation}

\subsection{Green's matrix $g(r, s)$ and explicit kernels}

Following \cite{Mou}, the Green's matrix $g(r, s)$ is given by
\begin{equation}\label{eq:formula_g}
	g(r, s) = \begin{cases}
		-\Phi(r) M_F \Phi^{-1}(s) + \Phi(r)\Phi^{-1}(s) = \Phi(r) M_E \Phi^{-1}(s), & A \le s \le r \le B, \\[0.25em]
		-\Phi(r) M_F \Phi^{-1}(s), & A \le r < s \le B.
	\end{cases}
\end{equation}
To express $g(r,s)$ compactly, define the auxiliary functions
\begin{equation}\label{eq:aux_funcs}
	\rho_1(x) := B + 1 - x, \quad \rho_\beta(x):= B + \beta - x, \quad \psi_1(x):= 1 - C_1(\alpha, \beta)(x - A), \quad \psi_2(x) := \beta - C_2(\alpha, \beta)(x - A).
\end{equation}
Then, $g(r, s)$ takes the following form:
\begin{enumerate}[$(i)$]
	\item For $r < s$ ($A \le r < s \le B$),
	\begin{equation}\label{eq:g_right}
		g(r, s) = \frac{1}{D(\alpha, \beta)} \begin{pmatrix}
			-\psi_1(r) &\quad  -\rho_1(s) \psi_1(r) &\quad  \alpha\beta \rho_1(r) &\quad  \alpha\beta \rho_1(r) \rho_\beta(s) \\[0.25em]
			C_1 &\quad  C_1 \rho_1(s) &\quad  -\alpha\beta &\quad  -\alpha\beta \rho_\beta(s) \\[0.25em]
			\alpha \rho_\beta(r) &\quad  \alpha \rho_\beta(r) \rho_1(s) &\quad  -\psi_2(r) &\quad  -\rho_\beta(s) \psi_2(r) \\[0.25em]
			-\alpha &\quad  -\alpha \rho_1(s) &\quad  C_2 &\quad  C_2 \rho_\beta(s)
		\end{pmatrix}.
	\end{equation}
	
	\item For $s \le r$ ($A \le s \le r \le B$),
	\begin{equation}\label{eq:g_left}
		g(r, s) = \frac{1}{D(\alpha, \beta)} \begin{pmatrix}
			C_1 \rho_1(r) & \quad  -\rho_1(r) \psi_1(s) & \quad  \alpha\beta \rho_1(r) & \quad   \alpha\beta \rho_1(r) \rho_\beta(s) \\[0.25em]
			C_1 & \quad \psi_1(s) & \quad -\alpha\beta & \quad -\alpha\beta \rho_\beta(s) \\[0.25em]
			\alpha \rho_\beta(r) & \quad \alpha \rho_\beta(r) \rho_1(s) &\quad  C_2 \rho_\beta(r) & \quad  -\rho_\beta(r) \psi_2(s) \\[0.25em]
			-\alpha & \quad -\alpha \rho_1(s) & \quad C_2 & \quad \psi_2(s)
		\end{pmatrix}.
	\end{equation}
\end{enumerate}
Note that on the diagonal $r = s$, the jump satisfies $g(s^+, s) - g(s^-, s) = \Phi(s)[M_E + M_F]\Phi^{-1}(s) = I_4$.

Finally, by \cite[Corollary 3.3]{Mou} the four scalar Green's kernels are identified as
\begin{equation}\label{eq:kernel_ident}
	G_1^1(r, s) := -g_{1,2}(r, s), \quad G_2^1(r, s) := -g_{1,4}(r, s), \quad G_1^2(r, s) :=-g_{3,2}(r, s), \quad G_2^2(r, s) := -g_{3,4}(r, s).
\end{equation}
Hence, evaluating these entries from \eqref{eq:g_right} and \eqref{eq:g_left} yields
\begin{align}
G_1^1(r, s) &= \frac{(B + 1 - \max\{r, s\}) [1 - C_1(\alpha, \beta)(\min\{r, s\} - A)]}{D(\alpha, \beta)}, \quad &&(r, s) \in I \times I, \label{eq:G11} \\[0.4em]
G_2^1(r, s) &= \frac{-\alpha\beta (B + 1 - r)(B + \beta - s)}{D(\alpha, \beta)}, \quad &&(r, s) \in I \times I, \label{eq:G21} \\[0.4em]
G_1^2(r, s) &= \frac{-\alpha (B + \beta - r)(B + 1 - s)}{D(\alpha, \beta)}, \quad && (r, s) \in I \times I, \label{eq:G12} \\[0.4em]
G_2^2(r, s) &= \frac{(B + \beta - \max\{r, s\}) [\beta - C_2(\alpha, \beta)(\min\{r, s\} - A)]}{D(\alpha, \beta)}, \quad &&(r, s) \in I \times I. \label{eq:G22}
\end{align}

\section{Monotonicity analysis and extrema of Green's functions in Example \ref{ex4.1}}
\label{app:B}
Let $A < B$ with unit domain length $L := B - A = 1$. We introduce the normalized spatial coordinates on the unit interval:
\begin{equation}\label{eq:coord_norm}
	\xi := r - A \in [0, 1], \qquad \eta := s - A \in [0, 1].
\end{equation}
Under these transformations the Green's kernels take the form
\begin{align}
	G_1^1(\xi, \eta) &= \frac{\bigl(2 - \max\{\xi, \eta\}\bigr) \bigl[1 - C_1(\alpha,\beta) \min\{\xi, \eta\}\bigr]}{D(\alpha,\beta)}, \label{eq:G11} \\[0.8em]
	G_2^1(\xi, \eta) &= \frac{-\alpha\beta \, (2 - \xi)(1 + \beta - \eta)}{D(\alpha,\beta)}, \label{eq:G12} \\[0.8em]
	G_1^2(\xi, \eta) &= \frac{-\alpha \, (1 + \beta - \xi)(2 - \eta)}{D(\alpha,\beta)}, \label{eq:G21} \\[0.8em]
	G_2^2(\xi, \eta) &= \frac{\bigl(1 + \beta - \max\{\xi, \eta\}\bigr) \bigl[\beta - C_2(\alpha,\beta) \min\{\xi, \eta\}\bigr]}{D(\alpha,\beta)}, \label{eq:G22}
\end{align}
for $(\xi, \eta) \in [0, 1] \times [0, 1]$, where $C_1\equiv C_1(\alpha,\beta) := 1 - (1 + \beta)\alpha^2,\, C_2\equiv C_2(\alpha,\beta) := 1 - 2\alpha^2$ and $	D\equiv D(\alpha,\beta) := \beta - (1 + \beta)C_2(\alpha, \beta)$.
Throughout this work, we assume
\begin{equation}\label{eq:hyp_standing}
	\beta \ge 0, \qquad \alpha < 0, \qquad \text{and} \qquad D(\alpha,\beta) > 0.
\end{equation}
Condition $D(\alpha,\beta) > 0$ is equivalent to $C_1(\alpha,\beta) < 1/2$, and if $\beta > 0$, it is also equivalent to $C_2(\alpha,\beta) < \frac{\beta}{1+\beta}$.

\subsection{Cross kernels \texorpdfstring{$G_1^2$}{G12} and \texorpdfstring{$G_2^1$}{G21}}

Under hypothesis \eqref{eq:hyp_standing}, $-\alpha > 0$, $-\alpha\beta \ge 0$, and the off-diagonal factors $2 - \xi \ge 1 > 0$, $1 + \beta - \xi \ge \beta \ge 0$, $2 - \eta \ge 1 > 0$, and $1 + \beta - \eta \ge \beta \ge 0$ are non-negative for all $(\xi, \eta) \in [0, 1]^2$. Hence, $G_1^2$ and $G_2^1$ are everywhere non-negative.

\begin{proposition}\label{prop:cross_kernels}
	Assume \eqref{eq:hyp_standing}. The off-diagonal kernels $G_1^2(\xi, \eta)$ and $G_2^1(\xi, \eta)$ are non-increasing with respect to each variable on the unit square $[0, 1] \times [0, 1]$. Consequently, their global extrema on $[0, 1] \times [0, 1]$ are given by
	\begin{align}
		M_2^1 &:= \max_{(\xi, \eta) \in [0, 1]^2} G_1^2(\xi, \eta) = G_1^2(0, 0) = \frac{-2\alpha\beta(1 + \beta)}{D(\alpha,\beta)}, \label{eq:M12} \\[0.4em]
		M_1^2 &:= \max_{(\xi, \eta) \in [0, 1]^2} G_2^1(\xi, \eta) = G_2^1(0, 0) = \frac{-2\alpha(1 + \beta)}{D(\alpha,\beta)}, \label{eq:M21} \\[0.4em]
		m_2^1 &:= \min_{(\xi, \eta) \in [0, 1]^2} G_1^2(\xi, \eta) = G_1^2(1, 1) = \frac{-\alpha\beta^2}{D(\alpha,\beta)}, \label{eq:m12} \\[0.4em]
		m_1^2 &:= \min_{(\xi, \eta) \in [0, 1]^2} G_2^1(\xi, \eta) = G_2^1(1, 1) = \frac{-\alpha\beta}{D(\alpha,\beta)}. \label{eq:m21}
	\end{align}
	Moreover, the optimal envelope functions $\Phi_i^j(s) := \max_{r \in [A, B]} G_i^j(r, s)$ are attained at $r = A$ ($\xi = 0$):
	\begin{align}
		\Phi_2^1(s) &= G_2^1(A, s) = \frac{-2\alpha\beta(B + \beta - s)}{D(\alpha, \beta)} = \frac{-2\alpha\beta(1 + \beta - \eta)}{D(\alpha, \beta)}, \label{eq:Phi12} \\[0.6em]
		\Phi_1^2(s) &= G_2^1(A, s) = \frac{-\alpha(1 + \beta)(B + 1 - s)}{D(\alpha, \beta)} = \frac{-\alpha(1 + \beta)(2 - \eta)}{D(\alpha, \beta)}. \label{eq:Phi21}
	\end{align}
\end{proposition}

\subsection{Primary kernel \texorpdfstring{$G_1^1$}{G11}}

The kernel $G_1^1(\xi, \eta)$ is symmetric (i.e. $G_1^1(\xi, \eta) = G_1^1(\eta, \xi)$). Its spatial derivative with respect to $\xi$ on each side of the diagonal $\xi = \eta$ is:
\begin{equation}\label{eq:dG11_dxi}
	\frac{\partial}{\partial \xi} G_1^1(\xi, \eta) = 
	\begin{cases}
		\dfrac{-C_1(\alpha,\beta)(2 - \eta)}{D(\alpha,\beta)}, & \text{if } \xi < \eta, \\[0.8em]
		\dfrac{-(1 - C_1(\alpha,\beta)\eta)}{D(\alpha,\beta)}, & \text{if } \xi > \eta.
	\end{cases}
\end{equation}
For any admissible $C_1(\alpha,\beta) < 1/2$, we have $1 - C_1(\alpha,\beta)\eta \ge 1 - \max(0, C_1) > 0$. Therefore, $\frac{\partial}{\partial \xi} G_1^1(\xi, \eta) < 0$ strictly for all $\xi > \eta$. On the other hand, the sign of the derivative for $\xi < \eta$ is governed entirely by $C_1(\alpha, \beta)$.

\begin{proposition}[Monotonicity of $G_1^1$]\label{prop:G11_monotonia}
	Assume \eqref{eq:hyp_standing}.
	\begin{enumerate}
		\item If $0 \le C_1(\alpha,\beta) < 1/2$, $G_1^1(\xi, \eta)$ is non-increasing in each variable on $[0, 1] \times [0, 1]$. Its global extrema are
		\begin{equation}\label{eq:extrema_G11_pos}
			M_1^1 := G_1^1(0, 0) = \frac{2}{D(\alpha,\beta)}, \qquad m_1^1 := G_1^1(1, 1) = \frac{1 - C_1(\alpha,\beta)}{D(\alpha,\beta)} = \frac{(1+\beta)\alpha^2}{D(\alpha,\beta)}.
		\end{equation}
		The optimal envelope is strictly decreasing:
		\begin{equation}\label{eq:Phi11_pos}
			\Phi_1^1(s) := \max_{r \in [A, B]} G_1^1(r, s) = G_1^1(A, s) = \frac{B + 1 - s}{D(\alpha, \beta)} = \frac{2 - \eta}{D(\alpha, \beta)}.
		\end{equation}
		
		\item If $C_1(\alpha,\beta) < 0$, for each fixed $\eta \in [0, 1]$, $\xi \mapsto G_1^1(\xi, \eta)$ is strictly increasing on $[0, \eta)$ and strictly decreasing on $(\eta, 1]$. The optimal envelope is attained on the diagonal:
		\begin{align}\label{eq:Phi11_neg}
			\Phi_1^1(s) &:= \max_{r \in [A, B]} G_1^1(r, s) = G_1^1(s, s) \nonumber \\
			&= \frac{(B + 1 - s)[1 - C_1(\alpha,\beta)(s - A)]}{D(\alpha, \beta)} = \frac{(2 - \eta)(1 - C_1 \eta)}{D(\alpha, \beta)}.
		\end{align}
	\end{enumerate}
\end{proposition}

\begin{proposition}[Global Extrema of $G_1^1$ when $C_1(\alpha,\beta) < 0$]\label{prop:extrema_G11_neg}
	Let $C_1(\alpha,\beta) < 0$.
	\begin{enumerate}
		\item The absolute global minimum is
		\begin{equation}\label{eq:m11_neg}
			m_1^1 := \min_{(\xi, \eta) \in [0, 1]^2} G_1^1(\xi, \eta) = \frac{1}{D(\alpha,\beta)},
		\end{equation}
		attained exclusively at the off-diagonal vertices $(\xi, \eta) = (0, 1)$ and $(\xi, \eta) = (1, 0)$.
		
		\item The absolute global maximum is attained on the diagonal $\xi = \eta$ and is given by
		\begin{equation}\label{eq:M11_neg}
			M_1^1 := \max_{(\xi, \eta) \in [0, 1]^2} G_1^1(\xi, \eta) = 
			\begin{cases}
				\dfrac{2}{D(\alpha,\beta)}, & \text{if } -\dfrac{1}{2} \le C_1(\alpha,\beta) < 0, \\[0.8em]
				\dfrac{D(\alpha,\beta)}{4|C_1(\alpha,\beta)|}, & \text{if } C_1(\alpha,\beta) < -\dfrac{1}{2}.
			\end{cases}
		\end{equation}
		It is attained at $(0, 0)$ if $-1/2 \le C_1 < 0$, and at $(\xi^*, \xi^*)$ with $\xi^* := 1 + \frac{1}{2C_1(\alpha,\beta)} \in (0, 1)$ if $C_1 < -1/2$.
	\end{enumerate}
\end{proposition}

\noindent\textit{Proof.}
	For $(i)$, since $G_1^1$ increases for $\xi \le \eta$ and decreases for $\xi \ge \eta$, the minimum over $\xi \in [0, 1]$ must occur at $\xi \in \{0, 1\}$. Evaluating $g_0(\eta) := G_1^1(0, \eta) = \frac{2 - \eta}{D}$ and $g_1(\eta) := G_1^1(1, \eta) = \frac{1 - C_1 \eta}{D}$, we observe that $g_0(\eta)$ is strictly decreasing in $\eta$ (minimum $g_0(1) = 1/D$), while $g_1(\eta)$ is strictly increasing in $\eta$ because $-C_1 > 0$ (minimum $g_1(0) = 1/D$). At the outer vertex $(\xi, \eta) = (1, 1)$, $G_1^1(1, 1) = \frac{1 - C_1}{D} > \frac{1}{D}$. Thus, $m_1^1 = 1/D$ uniquely at $(0, 1)$ and $(1, 0)$.
	
	For $(ii)$, the global maximum lies on the diagonal $f(\xi) := G_1^1(\xi, \xi) = \frac{(2 - \xi)(1 - C_1 \xi)}{D}$. We compute $f'(\xi) = \frac{-(1 + 2C_1) + 2C_1 \xi}{D}$ and $f''(\xi) = \frac{2C_1}{D} < 0$. If $-1/2 \le C_1 < 0$, then $-(1 + 2C_1) \le 0$ and $2C_1 \xi \le 0$, so $f'(\xi) \le 0$ on $[0, 1]$, placing the maximum at $\xi = 0$. If $C_1 < -1/2$, then $f'(0) > 0$ and $f'(1) = -1/D < 0$, yielding a unique interior critical point $\xi^* = 1 + \frac{1}{2C_1} \in (0, 1)$ with $f(\xi^*) = \frac{D}{4|C_1|}$.
\qed

\subsection{Secondary kernel \texorpdfstring{$G_2^2$}{G22}}

The kernel $G_2^2(\xi, \eta)$ is likewise symmetric. Its derivative with respect to $\xi$ is:
\begin{equation}\label{eq:dG22_dxi}
	\frac{\partial}{\partial \xi} G_2^2(\xi, \eta) = 
	\begin{cases}
		\dfrac{-C_2(\alpha,\beta)(1 + \beta - \eta)}{D(\alpha,\beta)}, & \text{if } \xi < \eta, \\[0.8em]
		\dfrac{-(\beta - C_2(\alpha,\beta)\eta)}{D(\alpha,\beta)}, & \text{if } \xi > \eta.
	\end{cases}
\end{equation}
When $\beta > 0$ and $C_2 \ge 0$, the condition $D > 0$ implies $\beta > \frac{C_2}{1 - C_2} \ge C_2 \ge 0$. Hence $\beta - C_2 \eta \ge \beta - C_2 > 0$. When $C_2 < 0$, $-C_2 > 0$, so $\beta - C_2 \eta \ge \beta > 0$. Thus, for $\xi > \eta$, $\frac{\partial}{\partial \xi} G_2^2 < 0$ strictly in all cases.

\begin{proposition}[Monotonicity of $G_2^2$ for $\beta > 0$]\label{prop:G22_monotonia}
	Assume \eqref{eq:hyp_standing} with $\beta > 0$.
	\begin{enumerate}
		\item If $C_2(\alpha,\beta) \ge 0$, $G_2^2(\xi, \eta)$ is non-increasing in each variable on $[0, 1] \times [0, 1]$. Its global extrema are
		\begin{equation}\label{eq:extrema_G22_pos}
			M_2^2 := G_2^2(0, 0) = \frac{\beta(1+\beta)}{D(\alpha,\beta)}, \qquad m_2^2 := G_2^2(1, 1) = \frac{\beta(\beta - C_2(\alpha,\beta))}{D(\alpha,\beta)}.
		\end{equation}
		The optimal envelope is affine and strictly decreasing:
		\begin{equation}\label{eq:Phi22_pos}
			\Phi_2^2(s) := \max_{r \in [A, B]} G_2^2(r, s) = G_2^2(A, s) = \frac{\beta(B + \beta - s)}{D(\alpha, \beta)} = \frac{\beta(1 + \beta - \eta)}{D(\alpha, \beta)}.
		\end{equation}
		
		\item If $C_2(\alpha,\beta) < 0$, for each fixed $\eta \in [0, 1]$, $\xi \mapsto G_2^2(\xi, \eta)$ is strictly increasing on $[0, \eta)$ and strictly decreasing on $(\eta, 1]$. The optimal envelope is attained on the diagonal:
		\begin{align}\label{eq:Phi22_neg}
			\Phi_2^2(s) &:= \max_{r \in [A, B]} G_2^2(r, s) = G_2^2(s, s) \nonumber \\
			&= \frac{(B + \beta - s)[\beta - C_2(\alpha,\beta)(s - A)]}{D(\alpha, \beta)} = \frac{(1 + \beta - \eta)(\beta - C_2 \eta)}{D(\alpha, \beta)}.
		\end{align}
	\end{enumerate}
\end{proposition}

\begin{proposition}[Global Extrema of $G_2^2$ when $\beta > 0$ and $C_2(\alpha,\beta) < 0$]\label{prop:extrema_G22_neg}
	Let $\beta > 0$ and $C_2(\alpha,\beta) < 0$.
	\begin{enumerate}
		\item The absolute global minimum is
		\begin{equation}\label{eq:m22_neg}
			m_2^2 := \min_{(\xi, \eta) \in [0, 1]^2} G_2^2(\xi, \eta) = \frac{\beta^2}{D(\alpha,\beta)},
		\end{equation}
		attained exclusively at the off-diagonal vertices $(\xi, \eta) = (0, 1)$ and $(\xi, \eta) = (1, 0)$.
		
		\item The absolute global maximum is attained on the diagonal $\xi = \eta$ and is given by
		\begin{equation}\label{eq:M22_neg}
			M_2^2 := \max_{(\xi, \eta) \in [0, 1]^2} G_2^2(\xi, \eta) = 
			\begin{cases}
				\dfrac{\beta(1+\beta)}{D(\alpha,\beta)}, & \text{if } \beta (1 - |C_2(\alpha,\beta)|) \ge|C_2(\alpha,\beta)|, \\[0.8em]
				\dfrac{\bigl[|C_2(\alpha,\beta)|(1 + \beta) + \beta\bigr]^2}{4|C_2(\alpha,\beta)| D(\alpha,\beta)}, & \text{if } \beta (1 - |C_2(\alpha,\beta)|)< |C_2(\alpha,\beta)|.
			\end{cases}
		\end{equation}
		It is attained at $(0, 0)$ if $\beta \ge \frac{|C_2|}{1 - |C_2|}$, and at $(\xi^*, \xi^*)$ with $\xi^* := \frac{1+\beta}{2} - \frac{\beta}{2|C_2(\alpha,\beta)|} \in (0, 1)$ if $\beta < \frac{|C_2|}{1 - |C_2|}$.
	\end{enumerate}
\end{proposition}

\noindent\textit{Proof.}
	For $(i)$, the minimum with respect to $\xi$ occurs at $\xi \in \{0, 1\}$. Evaluating $g_0(\eta) := G_2^2(0, \eta) = \frac{\beta(1 + \beta - \eta)}{D}$ (strictly decreasing, minimum $g_0(1) = \beta^2/D$) and $g_1(\eta) := G_2^2(1, \eta) = \frac{\beta(\beta - C_2 \eta)}{D}$ (strictly increasing since $-C_2 > 0$, minimum $g_1(0) = \beta^2/D$). At $(\xi, \eta) = (1, 1)$, $G_2^2(1, 1) = \frac{\beta(\beta - C_2)}{D} > \frac{\beta^2}{D}$. Thus, the minimum value is $\beta^2/D$ at $(0, 1)$ and $(1, 0)$.
	
	For $(ii)$, on the diagonal $f(\xi) := G_2^2(\xi, \xi) = \frac{(1 + \beta - \xi)(\beta - C_2 \xi)}{D}$, we have $f''(\xi) = \frac{2C_2}{D} < 0$. The vertex of the parabola is located at $\xi^* = \frac{1+\beta}{2} - \frac{\beta}{2|C_2|}$. Setting $\xi^* \le 0$ yields $(1+\beta)|C_2| \le \beta \iff \beta (1-|C_2|)\ge |C_2|$, in which case $f$ is decreasing on $[0, 1]$ and the maximum is $f(0) = \frac{\beta(1+\beta)}{D}$. When $\beta (1-|C_2|) < |C_2|$, $\xi^* \in (0, 1)$, and direct substitution into $f(\xi^*)$ yields the quadratic maximum in \eqref{eq:M22_neg}.
\qed

\section{Computation of constants \texorpdfstring{$c_i$ and $\tilde{c}_i$}{ci and c-tilde-i} in Example \ref{ex4.1}}
\label{app:C}

Let $[a, b] \subset [A, B]$ be a designated subinterval with $A \le a < b \le B$. We define the normalized bounds
\begin{equation}\label{eq:norm_subinterval}
	\hat{a} := a - A \in [0, 1], \qquad \hat{b} := b - A \in [0, 1], \qquad \text{with } 0 \le \hat{a} < \hat{b} \le 1.
\end{equation}
In spatial coordinates: $B + 1 - b = 2 - \hat{b}$, $B + \beta - b = 1 + \beta - \hat{b}$, and $a - A = \hat{a}$.

For each $i \in \{1, 2\}$, the global and localized Harnack-type cone constants are defined by
\begin{equation}\label{eq:Harnack_constants_def}
	\tilde{c}_i := \min_{j \in \{1,2\}} \inf_{(r,s) \in [A, B]^2} \frac{G_j^i(r, s)}{\Phi_j^i(s)}, \qquad
	c_i := \min_{j \in \{1,2\}} \inf_{\substack{r \in [a,b] \\ s \in [A, B]}} \frac{G_j^i(r, s)}{\Phi_j^i(s)}.
\end{equation}
In normalized coordinates $(\xi, \eta) \in [0, 1]^2$, these correspond to
\begin{equation}\label{eq:Harnack_norm_def}
	\tilde{c}_i = \min_{j \in \{1,2\}} \inf_{(\xi,\eta) \in [0, 1]^2} \frac{G_j^i(\xi, \eta)}{\Phi_j^i(\eta)}, \qquad
	c_i = \min_{j \in \{1,2\}} \inf_{\substack{\xi \in [\hat{a},\hat{b}] \\ \eta \in [0, 1]}} \frac{G_j^i(\xi, \eta)}{\Phi_j^i(\eta)}.
\end{equation}

We first evaluate the off-diagonal Green's kernels $G_1^2$ and $G_2^1$.

\begin{proposition}\label{prop:cross_ratios}
	Assume \eqref{eq:hyp_standing}.
	\begin{enumerate}
		\item For $i = 1, j = 2$, assuming $\beta > 0$:
		\begin{equation}\label{eq:ratio_G12}
			\frac{G_2^1(\xi, \eta)}{\Phi_2^1(\eta)} = \frac{2 - \xi}{2} = 1 - \frac{\xi}{2}.
		\end{equation}
		Consequently,
		\begin{equation}\label{eq:inf_G12}
			\inf_{(\xi, \eta) \in [0, 1]^2} \frac{G_2^1(\xi, \eta)}{\Phi_2^1(\eta)} = \frac{1}{2}, \qquad
			\inf_{\substack{\xi \in [\hat{a}, \hat{b}] \\ \eta \in [0, 1]}} \frac{G_2^1(\xi, \eta)}{\Phi_2^1(\eta)} = \frac{2 - \hat{b}}{2} = 1 - \frac{\hat{b}}{2} = \frac{B + 1 - b}{2}.
		\end{equation}
		
		\item For $i = 2, j = 1$, for all $\beta \ge 0$:
		\begin{equation}\label{eq:ratio_G21}
			\frac{G_1^2(\xi, \eta)}{\Phi_1^2(\eta)} = \frac{1 + \beta - \xi}{1 + \beta} = 1 - \frac{\xi}{1 + \beta}.
		\end{equation}
		Consequently,
		\begin{equation}\label{eq:inf_G21}
			\inf_{(\xi, \eta) \in [0, 1]^2} \frac{G_1^2(\xi, \eta)}{\Phi_1^2(\eta)} = \frac{\beta}{1 + \beta}, \qquad
			\inf_{\substack{\xi \in [\hat{a}, \hat{b}] \\ \eta \in [0, 1]}} \frac{G_1^2(\xi, \eta)}{\Phi_1^2(\eta)} = \frac{1 + \beta - \hat{b}}{1 + \beta} = 1 - \frac{\hat{b}}{1 + \beta} = \frac{B + \beta - b}{1 + \beta}.
		\end{equation}
	\end{enumerate}
\end{proposition}

\subsection{Derivation of the first constants \texorpdfstring{$c_1$}{c1} and \texorpdfstring{$\tilde{c}_1$}{c1 tilde}}

We now evaluate the ratio $L_1^1(\xi, \eta) := \frac{G_1^1(\xi, \eta)}{\Phi_1^1(\eta)}$.

\begin{lemma}[Ratio $G_1^1 / \Phi_1^1$ for $0 \le C_1(\alpha, \beta) < 1/2$]\label{lem:ratio_G11_pos}
	Assume $0 \le C_1(\alpha, \beta) < 1/2$. Then
	\begin{equation}\label{eq:ratio_G11_pos_form}
		\frac{G_1^1(\xi, \eta)}{\Phi_1^1(\eta)} = \begin{cases}
			1 - C_1(\alpha, \beta)\xi, & \text{if } \xi \le \eta, \\[0.6em]
			\dfrac{2 - \xi}{2 - \eta}\bigl(1 - C_1(\alpha, \beta)\eta\bigr), & \text{if } \xi \ge \eta.
		\end{cases}
	\end{equation}
	Moreover,
	\begin{equation}\label{eq:inf_G11_pos}
		\inf_{(\xi, \eta) \in [0, 1]^2} \frac{G_1^1(\xi, \eta)}{\Phi_1^1(\eta)} = \frac{1}{2}, \qquad
		\inf_{\substack{\xi \in [\hat{a}, \hat{b}] \\ \eta \in [0, 1]}} \frac{G_1^1(\xi, \eta)}{\Phi_1^1(\eta)} = \frac{2 - \hat{b}}{2} = \frac{B + 1 - b}{2}.
	\end{equation}
\end{lemma}

\noindent\textit{Proof.}
	Dividing \eqref{eq:G11} by \eqref{eq:Phi11_pos} gives \eqref{eq:ratio_G11_pos_form}. For any fixed $\eta \in [0, 1]$,  for $\xi \le \eta$, $\xi \mapsto 1 - C_1 \xi$ has derivative $-C_1 \le 0$, while for $\xi \ge \eta$, $\xi \mapsto \frac{2 - \xi}{2 - \eta}(1 - C_1 \eta)$ has derivative $-\frac{1 - C_1 \eta}{2 - \eta} < 0$.

	Thus, for every $\eta$, $\xi \mapsto L_1^1(\xi, \eta)$ is non-increasing across the entire interval $[0, 1]$.
	
	For the global infimum on $[0, 1]^2$, the minimum with respect to $\xi$ is attained at $\xi = 1$. Set
	\[
	L_1^1(1, \eta) = \frac{1 - C_1 \eta}{2 - \eta} =: h(\eta).
	\]
	Differentiating $h(\eta)$ on $[0, 1]$, we obtain
	\[
	h'(\eta) = \frac{-C_1(2 - \eta) - (1 - C_1 \eta)(-1)}{(2 - \eta)^2} = \frac{1 - 2C_1}{(2 - \eta)^2} = \frac{D(\alpha, \beta)}{(2 - \eta)^2} > 0.
	\]
	Thus $h(\eta)$ is strictly increasing in $\eta$. Its minimum on $[0, 1]$ is attained at $\eta = 0$. Hence,
	\[
	\inf_{(\xi, \eta) \in [0, 1]^2} L_1^1(\xi, \eta) = h(0) = \frac{1}{2}.
	\]
	
	For the localized infimum on $[\hat{a}, \hat{b}] \times [0, 1]$, since $L_1^1(\xi, \eta)$ is non-increasing in $\xi$, its minimum on $[\hat{a}, \hat{b}]$ is attained at $\xi = \hat{b}$. For $\eta \in [0, \hat{b}]$ (where $\hat{b} \ge \eta$), $L_1^1(\hat{b}, \eta) = (2 - \hat{b})h(\eta)$, so since $h'(\eta) > 0$, the minimum on $[0, \hat{b}]$ occurs at $\eta = 0$, giving $ \frac{2 - \hat{b}}{2}$. On the other hand, for $\eta \in [\hat{b}, 1]$ (where $\hat{b} \le \eta$), $L_1^1(\hat{b}, \eta) = 1 - C_1 \hat{b}$ and since $C_1 < 1/2$, we have $1 - C_1 \hat{b} \ge 1 - \frac{\hat{b}}{2} = \frac{2 - \hat{b}}{2}$. Therefore, the infimum over all $\eta \in [0, 1]$ is $\frac{2 - \hat{b}}{2} = \frac{B + 1 - b}{2}$.
\qed

\begin{lemma}[Ratio $G_1^1 / \Phi_1^1$ for $C_1(\alpha, \beta) < 0$]\label{lem:ratio_G11_neg}
	Assume $C_1(\alpha, \beta) < 0$. Then
	\begin{equation}\label{eq:ratio_G11_neg_form}
		\frac{G_1^1(\xi, \eta)}{\Phi_1^1(\eta)} = \begin{cases}
			\dfrac{1 - C_1(\alpha, \beta)\xi}{1 - C_1(\alpha, \beta)\eta}, & \text{if } \xi \le \eta, \\[0.8em]
			\dfrac{2 - \xi}{2 - \eta}, & \text{if } \xi \ge \eta.
		\end{cases}
	\end{equation}
	Moreover,
	\begin{equation}\label{eq:inf_G11_neg}
		\inf_{(\xi, \eta) \in [0, 1]^2} \frac{G_1^1(\xi, \eta)}{\Phi_1^1(\eta)} = \min\left\{ \frac{1}{2}, \, \frac{1}{1 - C_1(\alpha, \beta)} \right\},
	\end{equation}
	and
	\begin{align}\label{eq:inf_G11_neg_sub}
		\inf_{\substack{\xi \in [\hat{a}, \hat{b}] \\ \eta \in [0, 1]}} \frac{G_1^1(\xi, \eta)}{\Phi_1^1(\eta)} &= \min\left\{ \frac{2 - \hat{b}}{2}, \, \frac{1 - C_1(\alpha, \beta)\hat{a}}{1 - C_1(\alpha, \beta)} \right\} \nonumber \\
		&= \min\left\{ \frac{B + 1 - b}{2}, \, \frac{1 - C_1(\alpha, \beta)(a - A)}{1 - C_1(\alpha, \beta)} \right\}.
	\end{align}
\end{lemma}

\noindent\textit{Proof.}
	Dividing \eqref{eq:G11} by \eqref{eq:Phi11_neg} yields \eqref{eq:ratio_G11_neg_form}. For each fixed $\eta \in [0, 1]$, For $\xi \le \eta$, $\xi \mapsto \frac{1 - C_1 \xi}{1 - C_1 \eta}$ is strictly increasing in $\xi$ because $-C_1 > 0$, while for $\xi \ge \eta$, $\xi \mapsto \frac{2 - \xi}{2 - \eta}$ is strictly decreasing in $\xi$. Thus, 	the maximum is attained at $\xi = \eta$ where the ratio equals 1 and, on any subinterval, the minimum in $\xi$ must be achieved at an endpoint.
	
	For the computation of the global infimum on $[0, 1]^2$, note that:
	\begin{enumerate}[$(i)$]
	\item For $\xi \le \eta$, the minimum with respect to $\xi$ is achieved at $\xi = 0$, yielding $\frac{1}{1 - C_1 \eta}$. Since $-C_1 > 0$, the denominator increases with $\eta$, so the function is decreasing in $\eta$; hence, its minimum on $\eta \in [0, 1]$ is attained at $\eta = 1$, which gives $\frac{1}{1 - C_1}$.
	\item For $\xi \ge \eta$, the minimum with respect to $\xi$ is achieved at $\xi = 1$, yielding $\frac{2 - 1}{2 - \eta} = \frac{1}{2 - \eta}$. Because the denominator $2 - \eta$ decreases with $\eta$, the function increases in $\eta$; thus, its minimum on $\eta \in [0, 1]$ occurs at $\eta = 0$, which gives $\frac{1}{2}$.
	\end{enumerate}
	
	For the localized infimum on $[\hat{a}, \hat{b}] \times [0, 1]$, for any fixed $\eta$, the minimum over $\xi \in [\hat{a}, \hat{b}]$ must occur at $\xi = \hat{a}$ or $\xi = \hat{b}$. Observe that:
	\begin{enumerate}[$(i)$]
		\item At $\eta = 0$, every $\xi \in [\hat{a}, \hat{b}]$ satisfies $\xi \ge \eta = 0$. Hence $L_1^1(\xi, 0) = \frac{2 - \xi}{2}$, which is strictly decreasing in $\xi$. Its minimum on $[\hat{a}, \hat{b}]$ occurs at $\xi = \hat{b}$, giving $\frac{2 - \hat{b}}{2}$.
		\item At $\eta = 1$, every $\xi \in [\hat{a}, \hat{b}]$ satisfies $\xi \le \eta = 1$. Hence $L_1^1(\xi, 1) = \frac{1 - C_1 \xi}{1 - C_1}$, which is strictly increasing in $\xi$. Its minimum on $[\hat{a}, \hat{b}]$ occurs at $\xi = \hat{a}$, giving $\frac{1 - C_1 \hat{a}}{1 - C_1}$.
		\item For any intermediate $\eta \in (0, 1)$, by monotonicity of each linear piece, the values $L_1^1(\hat{a}, \eta)$ and $L_1^1(\hat{b}, \eta)$ are bounded below by their values at the boundary $\eta \in \{0, 1\}$, meaning
		\[
		\min_{\xi \in [\hat{a}, \hat{b}]} L_1^1(\xi, \eta) \ge \min\left\{ L_1^1(\hat{b}, 0), \, L_1^1(\hat{a}, 1) \right\} = \min\left\{ \frac{2 - \hat{b}}{2}, \, \frac{1 - C_1 \hat{a}}{1 - C_1} \right\}.
		\]
	\end{enumerate}
	Thus, the infimum over all $\xi \in [\hat{a}, \hat{b}]$ and $\eta \in [0, 1]$ is precisely \eqref{eq:inf_G11_neg_sub}.
\qed

We are now ready to state the full characterization of the constants $\tilde{c}_1$ and $c_1$.

\begin{theorem}[Explicit Formulas for $\tilde{c}_1$ and $c_1$]\label{thm:c1_final}
	Assume \eqref{eq:hyp_standing} with $\beta > 0$. The Harnack constants for the first component are given by
		\begin{equation}\label{eq:c1_tilde_piecewise}
			\tilde{c}_1 = \min\left\{ \frac{1}{2}, \, \frac{1}{1 - C_1(\alpha, \beta)} \right\} = \begin{cases}
				\dfrac{1}{2}, & \text{if } -1 \le C_1(\alpha, \beta) < \dfrac{1}{2}, \\[0.8em]
				\dfrac{1}{1 - C_1(\alpha, \beta)} = \dfrac{1}{(1 + \beta)\alpha^2}, & \text{if } C_1(\alpha, \beta) < -1,
			\end{cases}
		\end{equation}
		and
		\begin{equation}\label{eq:c1_localized_piecewise}
			c_1 = \begin{cases}
				\dfrac{B + 1 - b}{2}, & \text{if } 0 \le C_1(\alpha, \beta) < \dfrac{1}{2}, \\[0.8em]
				\min\left\{ \dfrac{B + 1 - b}{2}, \, \dfrac{1 - C_1(\alpha, \beta)(a - A)}{1 - C_1(\alpha, \beta)} \right\}, & \text{if } C_1(\alpha, \beta) < 0.
			\end{cases}
		\end{equation}
	In particular, $\tilde{c}_1 > 0$ and $c_1 > 0$ for all admissible parameter values.
\end{theorem}

\subsection{Derivation of the second constants \texorpdfstring{$c_2$}{c2} and \texorpdfstring{$\tilde{c}_2$}{c2 tilde}}

Next, we evaluate the ratio $L_2^2(\xi, \eta) := \frac{G_2^2(\xi, \eta)}{\Phi_2^2(\eta)}$ for $\beta > 0$.

\begin{lemma}[Ratio $G_2^2 / \Phi_2^2$ for $\beta > 0$ and $C_2(\alpha, \beta) \ge 0$]\label{lem:ratio_G22_pos}
	Assume $\beta > 0$ and $C_2(\alpha, \beta) \ge 0$. Then
	\begin{equation}\label{eq:ratio_G22_pos_form}
		\frac{G_2^2(\xi, \eta)}{\Phi_2^2(\eta)} = \begin{cases}
			1 - \dfrac{C_2(\alpha, \beta)}{\beta}\xi, & \text{if } \xi \le \eta, \\[0.8em]
			\dfrac{1 + \beta - \xi}{1 + \beta - \eta}\left(1 - \dfrac{C_2(\alpha, \beta)}{\beta}\eta\right), & \text{if } \xi \ge \eta.
		\end{cases}
	\end{equation}
	Moreover,
	\begin{equation}\label{eq:inf_G22_pos}
		\inf_{(\xi, \eta) \in [0, 1]^2} \frac{G_2^2(\xi, \eta)}{\Phi_2^2(\eta)} = \frac{\beta}{1 + \beta}, \qquad
		\inf_{\substack{\xi \in [\hat{a}, \hat{b}] \\ \eta \in [0, 1]}} \frac{G_2^2(\xi, \eta)}{\Phi_2^2(\eta)} = \frac{1 + \beta - \hat{b}}{1 + \beta} = \frac{B + \beta - b}{1 + \beta}.
	\end{equation}
\end{lemma}

\noindent\textit{Proof.}
	Dividing \eqref{eq:G22} by \eqref{eq:Phi22_pos} yields \eqref{eq:ratio_G22_pos_form}. For any fixed $\eta \in [0, 1]$, for $\xi \le \eta$, $\xi \mapsto 1 - \frac{C_2}{\beta}\xi$ has derivative $-\frac{C_2}{\beta} \le 0$ while for $\xi \ge \eta$, $\xi \mapsto \frac{1 + \beta - \xi}{1 + \beta - \eta}(1 - \frac{C_2}{\beta}\eta)$ has derivative $-\frac{\beta - C_2 \eta}{\beta(1 + \beta - \eta)} < 0$. Thus, for every $\eta$, $\xi \mapsto L_2^2(\xi, \eta)$ is non-increasing across the entire unit interval $[0, 1]$.
	
	For the global infimum on $[0, 1]^2$, the minimum in $\xi$ occurs at $\xi = 1$. Set,
	\[
	L_2^2(1, \eta) = \frac{\beta}{1 + \beta - \eta} \frac{\beta - C_2 \eta}{\beta} = \frac{\beta - C_2 \eta}{1 + \beta - \eta} =: k(\eta).
	\]
	Differentiating $k(\eta)$ on $[0, 1]$,
	\[
	k'(\eta) = \frac{-C_2(1 + \beta - \eta) - (\beta - C_2 \eta)(-1)}{(1 + \beta - \eta)^2} = \frac{\beta - (1 + \beta)C_2}{(1 + \beta - \eta)^2} = \frac{D}{(1 + \beta - \eta)^2} > 0.
	\]
	Because $k'(\eta) > 0$, $k(\eta)$ is strictly increasing in $\eta$. Its minimum on $[0, 1]$ is achieved at $\eta = 0$. Hence,
	\[
	\inf_{(\xi, \eta) \in [0, 1]^2} L_2^2(\xi, \eta) = k(0) = \frac{\beta}{1 + \beta}.
	\]
	
	For the localized infimum on $[\hat{a}, \hat{b}] \times [0, 1]$, the minimum in $\xi$ is achieved at $\xi = \hat{b}$. For $\eta \in [0, \hat{b}]$, $L_2^2(\hat{b}, \eta) = \frac{1 + \beta - \hat{b}}{\beta} k(\eta)$ and since $k'(\eta) > 0$, the minimum on $[0, \hat{b}]$ occurs at $\eta = 0$, giving $\frac{1 + \beta - \hat{b}}{\beta} k(0) = \frac{1 + \beta - \hat{b}}{1 + \beta}$. On the other hand, for $\eta \in [\hat{b}, 1]$, $L_2^2(\hat{b}, \eta) = 1 - \frac{C_2}{\beta}\hat{b}$ and since $k(\eta)$ increases, then $1 - \frac{C_2}{\beta}\hat{b} = L_2^2(\hat{b}, \hat{b}) \ge L_2^2(\hat{b}, 0) = \frac{1 + \beta - \hat{b}}{1 + \beta}$. Thus, the infimum over all $\eta \in [0, 1]$ is $\frac{1 + \beta - \hat{b}}{1 + \beta} = \frac{B + \beta - b}{1 + \beta}$.
\qed

\begin{lemma}[Ratio $G_2^2 / \Phi_2^2$ for $\beta > 0$ and $C_2(\alpha, \beta) < 0$]\label{lem:ratio_G22_neg}
	Assume $\beta > 0$ and $C_2(\alpha, \beta) < 0$. Then
	\begin{equation}\label{eq:ratio_G22_neg_form}
		\frac{G_2^2(\xi, \eta)}{\Phi_2^2(\eta)} = \begin{cases}
			\dfrac{\beta - C_2(\alpha, \beta)\xi}{\beta - C_2(\alpha, \beta)\eta}, & \text{if } \xi \le \eta, \\[0.8em]
			\dfrac{1 + \beta - \xi}{1 + \beta - \eta}, & \text{if } \xi \ge \eta.
		\end{cases}
	\end{equation}
	Moreover,
	\begin{equation}\label{eq:inf_G22_neg}
		\inf_{(\xi, \eta) \in [0, 1]^2} \frac{G_2^2(\xi, \eta)}{\Phi_2^2(\eta)} = \min\left\{ \frac{\beta}{1 + \beta}, \, \frac{\beta}{\beta - C_2(\alpha, \beta)} \right\},
	\end{equation}
	and
	\begin{align}\label{eq:inf_G22_neg_sub}
		\inf_{\substack{\xi \in [\hat{a}, \hat{b}] \\ \eta \in [0, 1]}} \frac{G_2^2(\xi, \eta)}{\Phi_2^2(\eta)} &= \min\left\{ \frac{1 + \beta - \hat{b}}{1 + \beta}, \, \frac{\beta - C_2(\alpha, \beta)\hat{a}}{\beta - C_2(\alpha, \beta)} \right\} \nonumber \\
		&= \min\left\{ \frac{B + \beta - b}{1 + \beta}, \, \frac{\beta - C_2(\alpha, \beta)(a - A)}{\beta - C_2(\alpha, \beta)} \right\}.
	\end{align}
\end{lemma}

\noindent\text{Proof.}
	Dividing \eqref{eq:G22} by \eqref{eq:Phi22_neg} yields \eqref{eq:ratio_G22_neg_form}. For each fixed $\eta \in [0, 1]$, if $\xi \le \eta$, $\xi \mapsto \frac{\beta - C_2 \xi}{\beta - C_2 \eta}$ is strictly increasing in $\xi$ because $-C_2 > 0$ while if $\xi \ge \eta$, $\xi \mapsto \frac{1 + \beta - \xi}{1 + \beta - \eta}$ is strictly decreasing in $\xi$.
	
	For the computation of the global infimum on $[0, 1]^2$, note that:
	\begin{enumerate}[$(i)$]
		\item In $\xi \le \eta$, the minimum in $\xi$ is at $\xi = 0$, giving $\frac{\beta}{\beta - C_2 \eta}$. The denominator increases with $\eta$, so the minimum over $\eta \in [0, 1]$ occurs at $\eta = 1$, which gives $\frac{\beta}{\beta - C_2}$.
		\item In $\xi \ge \eta$, the minimum in $\xi$ is at $\xi = 1$, giving $\frac{1 + \beta - 1}{1 + \beta - \eta} = \frac{\beta}{1 + \beta - \eta}$. The denominator decreases with $\eta$, so the minimum over $\eta \in [0, 1]$ occurs at $\eta = 0$, which gives $\frac{\beta}{1 + \beta}$.
	\end{enumerate}
	Taking the minimum of both yields \eqref{eq:inf_G22_neg}.
	
	For the localized infimum on $[\hat{a}, \hat{b}] \times [0, 1]$, observe that:
	\begin{enumerate}[$(i)$]
		\item At $\eta = 0$, for all $\xi \ge 0$, $L_2^2(\xi, 0) = \frac{1 + \beta - \xi}{1 + \beta}$ is decreasing in $\xi$. Its minimum on $[\hat{a}, \hat{b}]$ is attained at $\xi = \hat{b}$, giving $\frac{1 + \beta - \hat{b}}{1 + \beta}$.
		\item At $\eta = 1$, for all $\xi \le 1$, $L_2^2(\xi, 1) = \frac{\beta - C_2 \xi}{\beta - C_2}$ is increasing in $\xi$. Its minimum on $[\hat{a}, \hat{b}]$ is at $\xi = \hat{a}$, giving $\frac{\beta - C_2 \hat{a}}{\beta - C_2}$.
		\item For any intermediate $\eta \in (0, 1)$, by monotonicity of each piece, the lower bound is again controlled by the boundaries, meaning
		\[
		\min_{\xi \in [\hat{a}, \hat{b}]} L_2^2(\xi, \eta) \ge \min\left\{ L_2^2(\hat{b}, 0), \, L_2^2(\hat{a}, 1) \right\} = \min\left\{ \frac{1 + \beta - \hat{b}}{1 + \beta}, \, \frac{\beta - C_2 \hat{a}}{\beta - C_2} \right\}.
		\]
	\end{enumerate}
	Hence, the infimum on $[\hat{a}, \hat{b}] \times [0, 1]$ is precisely \eqref{eq:inf_G22_neg_sub}.
\qed

We can now state the full characterization of the constants $\tilde{c}_2$ and $c_2$.

\begin{theorem}[Explicit Formulas for $\tilde{c}_2$ and $c_2$]\label{thm:c2_final}
	Assume \eqref{eq:hyp_standing} with $\beta > 0$. The Harnack constants for the second component are given by:
		\begin{equation}\label{eq:c2_tilde_piecewise}
			\tilde{c}_2 = \min\left\{ \frac{\beta}{1 + \beta}, \, \frac{\beta}{\beta - C_2(\alpha, \beta)} \right\} = \begin{cases}
				\dfrac{\beta}{1 + \beta}, & \text{if } -1 \le C_2(\alpha, \beta) < \dfrac{\beta}{1 + \beta}, \\[0.8em]
				\dfrac{\beta}{\beta - C_2(\alpha, \beta)}, & \text{if } C_2(\alpha, \beta) < -1.
			\end{cases}
		\end{equation}
		
		\begin{equation}\label{eq:c2_localized_piecewise}
			c_2 = \begin{cases}
				\dfrac{B + \beta - b}{1 + \beta}, & \text{if } C_2(\alpha, \beta) \ge 0, \\[0.8em]
				\min\left\{ \dfrac{B + \beta - b}{1 + \beta}, \, \dfrac{\beta - C_2(\alpha, \beta)(a - A)}{\beta - C_2(\alpha, \beta)} \right\}, & \text{if } C_2(\alpha, \beta) < 0.
			\end{cases}
		\end{equation}
	In particular, $\tilde{c}_2 > 0$ and $c_2 > 0$ for all admissible parameter values when $\beta > 0$.
\end{theorem}

\begin{remark}[The Dirichlet Boundary Case $\beta = 0$]\label{rem:dirichlet}
	When $\beta = 0$, the kernel 
	\[G_2^2(\xi, \eta) = (1 - \max(\xi, \eta))\min(\xi, \eta), \quad (\xi, \eta)\in[0,1]^2,\]
	 reduces to the classical Dirichlet Green's function, whose envelope is 	  $\Phi_2^2(s) = (1 - \eta)\eta$. The ratio is
	 \[L^2_2(\xi, \eta)=\begin{cases}
	 	\frac{\xi}{\eta}, \quad &\text{ for } \xi\leq \eta,\\
	 	\frac{1-\xi}{1-\eta},  \quad &\text{ for } \xi> \eta.
	 \end{cases}\]
 Because the Dirichlet kernel vanishes at the boundaries $\xi \in \{0, 1\}$, the global infimum on $[0, 1]^2$ is $\tilde{c}_2 = 0$. However, on any strictly interior subinterval $A < a < b < B$ ($0 < \hat{a} < \hat{b} < 1$), the localized constant remains strictly positive
	\begin{equation}
		c_2 = \min\{ \hat{a}, \, 1 - \hat{b} \} = \min\{ a - A, \, B - b \} > 0,
	\end{equation}
	enabling the application of Krasnosel'ski\u{\i}'s fixed point theorem in localized cones of the form \eqref{eq:cone_Ki}.
\end{remark}

\end{document}